\documentclass[12pt]{amsart}
\usepackage {amsfonts,amscd,amssymb,amsmath,tabularx,t1enc}
\usepackage{enumerate}

\begin{document}

\title{Ad-nilpotent ideals of a Parabolic subalgebra}
\author{C\'eline RIGHI}
\address{UMR 6086 CNRS, Département de Mathématiques, Téléport 2 - BP
  30179, Boulevard Marie et Pierre Curie, 86962 Futuroscope
  Chasseneuil Cedex, France}
\email{celine.righi@math.univ-poitiers.fr}

\begin{abstract}
We extend the results of Cellini-Papi (\cite{CP}, \cite{CP2}) on the 
characterizations of ad-nilpotent and abelian ideals of a Borel
subalgebra to parabolic subalgebras of a simple Lie algebra. These
characterizations are given in terms of elements of the affine Weyl
group and faces of alcoves. In the case of a parabolic subalgebra of a
classical simple Lie algebra, we give formulas for the number of these ideals.
\end{abstract}

%\ms{17B20}

\maketitle

\newtheorem {prop}{Proposition}[section]
\newtheorem {lemme}[prop]{Lemma}
\newtheorem {theoreme}[prop]{Theorem}
\newtheorem {corollaire}[prop]{Corollary}
\newtheorem {definition}[prop]{Definition}
\newtheorem {exemple}[prop]{Example}
\newtheorem {exemples}[prop]{Examples}
\newtheorem {remarque}[prop]{Remark}
\newtheorem {remarques}[prop]{Remarks}
\newtheorem {notations}[prop]{Notations}

\newcommand\lmaj{\mathbf{L}}
\newcommand\cset{\mathbb C}
\newcommand\zset{\mathbb Z}
\newcommand\nset{\mathbb N}
\newcommand\rset{\mathbb R}
\newcommand\racset {\widehat \Delta}
\newcommand\racbase {\widehat \Pi}
\newcommand\ilie {\mathfrak {i}}
\newcommand\blie {\mathfrak {b}}
\newcommand\glie {\mathfrak {g}}
\newcommand\hlie {\mathfrak {h}}
\newcommand\plie {\mathfrak {p}}
\newcommand\nlie {\mathfrak {n}}
\newcommand\calp {{\cal P}}
\newcommand \calf {{\mathcal F}}
\newcommand \calc {{\mathcal C}}
\newcommand \obarre {{\overline{\omega}}}
\newcommand \abelien {\mathcal{A}\mathfrak{b}}
\newcommand \aff {{\mathit {aff}}}
\newcommand \caln {{\mathcal N}}

\def\htrait{\omit\hrulefill}
\def\vtrait{\vrule height 10pt depth 6pt} 
\def\sboite#1{\vtrait\hbox to 16pt{\hfil #1\hfil}} 
\def\hboite#1{\hbox to 16pt{\hfil #1\hfil}} 

\def\vt{\vrule height 10pt depth 6pt} 
\def\sb#1{\vt\hbox to 36pt{\hfil #1\hfil}} 
\def\hb#1{\hbox to 36pt{\hfil #1\hfil}}

% The environment Dynkin
\newenvironment{Dynkin}{\setlength{\unitlength}{1.5pt}\begin{array}{l}}{\end{array}}

% The command \Dbloc{args}
\newcommand{\Dbloc}[1]{\begin{picture}(20,20)#1\end{picture}}
% Possible arguments inside \Dbloc
\newcommand{\Dcirc}{\put(10,10){\circle{4}}}
\newcommand{\Dbullet}{\put(10,10){\circle*{4}}}
\newcommand{\Deast}{\put(12,10){\line(1,0){8}}}
\newcommand{\Dwest}{\put(8,10){\line(-1,0){8}}}
\newcommand{\Dnorth}{\put(10,12){\line(0,1){8}}}
\newcommand{\Dsouth}{\put(10,8){\line(0,-1){8}}}
\newcommand{\Dnortheast}{\put(20,20){\line(-1,-1){8.6}}}
\newcommand{\Dnorthwest}{\put(0,20){\line(1,-1){8.6}}}
\newcommand{\Dsoutheast}{\put(20,0){\line(-1,1){8.6}}}
\newcommand{\Dsouthwest}{\put(0,0){\line(1,1){8.6}}}
\newcommand{\Ddoubleeast}{\put(10,12){\line(1,0){10}}\put(10,8){\line(1,0){10}}}
\newcommand{\Ddoublewest}{\put(10,12){\line(-1,0){10}}\put(10,8){\line(-1,0){10}}}
\newcommand{\Ddots}{\put(5,10){\circle*{0.9}}\put(10,10){\circle*{0.9}}\put(15,10){\circle*{0.9}}}
\newcommand{\Dtext}[2]{\makebox(20,20)[#1]{\scriptsize $#2$}}

% Other commands outside \Dbloc
\newcommand{\Dskip}{\\ [-4.5pt]}
\newcommand{\Dspace}{\Dbloc{}}
\newcommand{\Dleftarrow}{\hskip-5pt{\makebox(20,20)[l]{\Large$<$}}\hskip-25pt}
\newcommand{\Drightarrow}{\hskip-5pt{\makebox(20,20)[l]{\Large$>$}}\hskip-25pt}

%%% End of routines %%%

\section{Introduction}
Let $\glie$ be a complex simple Lie algebra of rank $l$. Let $\hlie$ be a  Cartan 
subalgebra and $\Delta$ the associated root system. We fix a system of
positive roots $\Delta^+$. Denote by $\Pi=\{\alpha_1,
\dots, \alpha_l\}$ the corresponding set of simple roots. Let $V$ be the Euclidian space $\sum_{k=1}^{l}\rset\alpha_k$. For each 
$\alpha \in \Delta$, let $ \mathfrak{g}_{ \alpha}$ 
be the root space of $\glie$ relative to $\alpha$.

For  $I \subset \Pi$, set $\Delta_I=\zset I \cap\Delta$.
We fix the corresponding standard parabolic subalgebra :   
$$
\plie_I=\hlie\oplus\left(\bigoplus\limits_{\alpha \in \Delta_I \cup \Delta^+}\glie_{\alpha}\right).
$$

An ideal $\ilie $ of $\plie_I$ is ad-nilpotent if and only if for all $x\in \ilie$, $ad_{\plie_I} x$ is nilpotent. Since any ideal of $\plie_I$ is $\hlie$-stable, we can deduce easily that an ideal is ad-nilpotent if and only if it is nilpotent. Moreover, we have $\ilie =\bigoplus\limits_{\alpha \in \Phi} \glie_{\alpha}$, for some subset $\Phi\subset \Delta^+ \setminus \Delta_I$. 

The purpose of this paper is to characterize and to enumerate 
ad-nilpotent and abelian ideals of a parabolic subalgebra.

When $I=\emptyset$,  $\plie_{\emptyset}$ is a Borel
subalgebra of $\glie$. Peterson proved that  the number of abelian
ideals of $\plie_{\emptyset}$ is $2^l$. Motived by this result, 
Cellini-Papi, Kostant, Panyushev, Sommers and Suter among others 
studied ad-nilpotent and abelian ideals of a Borel subalgebra.

In their articles \cite{CP} and \cite{CP2}, Cellini-Papi established
different characterizations of the set $\mathcal I$ of ad-nilpotent ideals of a Borel
subalgebra. They constructed a bijection between $\mathcal I$ and
certain elements of the affine Weyl group $\widehat W$ associated to $\Delta$,
which we shall call $\emptyset$-compatible. These 
$\emptyset$-compatible elements are in turn characterized by elements of the
coroot lattice. They established also, when $\glie$ is of classical
type, a correspondence between  ad-nilpotent ideals of $\glie$ and
some diagrams. We extend here their theory to the case of parabolic subalgebras.

Fix $I\subset \Pi$, we establish a bijection between
ad-nilpotent ideals of $\plie_{I}$ and what we call $I$-compatible elements of the
affine Weyl group $\widehat W$. We identify $\widehat W$ with the
group of affine transformations $W_{\aff}$ defined in \cite{Bourbaki} and we give a
characterization of the $I$-compatible elements via the dimension of
the intersection of the image of the fundamental alcove associated to
$W_{\aff}$ with some affine hyperplanes of $V$.

Using this result, we obtain an identity (theorem \ref{thabelien}) 
which generalizes the result of
Peterson. This identity links the
number of abelian ideals and the coefficients of the simple roots in
the highest root of $\Delta$. This allows us to conclude that if
$\glie$ is of type $A$ or $C$, the number of abelian ideals of
$\plie_I$ is $2^{l-\sharp I}$. It also explains why this result does
not hold in general.

On the other hand, the enumeration of ad-nilpotent and abelian ideals of 
$\plie_I$, when $\glie$ is
of classical type, is obtained using the diagrams given in \cite{CP},
modified, by deleting some rows and columns and grouping together some
boxes, according to the type of $\glie$. The formulas obtained depend
on the decomposition in connected components of $I$. Note
that the formulas obtained when $\glie$ is of type $A$ or $C$ are
nicer than the ones obtained when $\glie$ is of type $B$ or $D$
(theorems \ref{thabelien}, \ref{thenumeration} and propositions
\ref{propenumerationB}, \ref{propenumerationD}).

This paper is organized as follows : in section 2, we recall some
results on the affine Weyl group. In section 3, we give different
characterizations of $I$-compatible elements of $\widehat W$. 
The study, in section 4, 
of the volume of the
intersection of some affine hyperplanes on $V$ gives the 
results stated above on abelians ideals. Section 5 deals with the enumeration of
both ad-nilpotent and abelians ideals when $\glie$ is of classical type,
using diagrams. We give some remarks concerning the exceptional
cases and the relations with antichains in section 6.

{\bf Acknowledgment.} The author would like to thank the referee for her/his comments and remarks, and in particular for indicating the correct approach to the proof of Proposition \ref{alcove}.

%------------------------------------------------------------------------------------------------------------------------------------------------

\section{Generalities on the affine Weyl group}
We shall conserve the notations given in the introduction. In this
section, we shall recall some basic facts on the
affine Weyl group associated to $\Delta$. In particular, we need to
recall  two different realizations of this group. See \cite{Bourbaki}, 
\cite{CP} and \cite{Kac} for more
details.

We fix a scalar product $(.,.)$ on $V$. For $\alpha \in \Delta$, let 
$$
\alpha^{\vee}=\frac{2\alpha}{(\alpha,\alpha)}
$$
denote the corresponding
coroot. Denote by $ Q^{ \vee}$ the coroot 
lattice of $\Delta$. 

Let $W$ denote the Weyl group associated to $\Delta$. We shall realize the affine Weyl group as a group of automorphisms of the affine
root system associated to $\Delta$. Let $\widehat V = V \oplus \rset\delta \oplus
\rset\lambda$. We extend the above bilinear form on $V$ to a
non-degenerate symmetric bilinear form on $\widehat V$, also denoted $ (.,.)$, by
setting : 
$$
\ (\lambda,\lambda) = (\delta,\delta) = (\lambda,V)
=(\delta,V)= 0
\mbox{ and } (\delta,\lambda) = 1.
$$ 
Let $\widehat{\Delta}=\Delta
+\zset\delta$ be the set of (real) affine roots. We fix the
following positive root system 
$ \racset^+ = (\Delta^+ +\nset \delta) \cup(\Delta^- + \nset^{\ast}
\delta)$. We shall write $\alpha>0$ (resp. $\alpha<0$) if $\alpha \in
\widehat{\Delta}^+$ (resp. if $\alpha \in
\widehat{\Delta}^-=-\widehat{\Delta}^+$). Let $\theta$ be the highest root of 
$\Delta$, then $\racbase = \{\alpha_0=-\theta +\delta, \alpha_1,
\cdots, \alpha_l\}$ 
is the 
set of simple roots for $\widehat{\Delta}^+$.

Note that for any element $\beta+k\delta \in
\widehat{\Delta}^+$, we have
$(\beta+k\delta,\beta+k\delta)=(\beta,\beta)\not=0$. 
For all $\alpha \in \racset^+$,  we denote by $s_{\alpha}$ the reflection of
$\widehat V$ defined by
$$
s_{\alpha}(x) = x-\frac{2(\alpha,x)}{(\alpha,\alpha)}\alpha
$$ 
for $x \in \widehat V$.
The affine Weyl group $\widehat W$ is the subgroup of $\mathrm{Aut}(\widehat
V)$ generated by
$ \{ s_{\alpha};\ \alpha \in \racbase\}$. Observe that $w(\delta) = \delta
\ \mbox {for all} \ w\in \widehat W $, $s_{\alpha}(\lambda)=\lambda$,
 for all $\alpha \in \Pi$ and
 $s_{\alpha_0}(\lambda)=\lambda-\frac{2}{\|\theta\|^2}\alpha_0$, where
 $\|\theta\|=\sqrt{(\theta,\theta)}$.

Let $\tau \in Q^{\vee}$, we define the endomorphism $t_{\tau}$ of $\widehat V$ by :

\begin{eqnarray}\label{deftranslation}
t_{\tau}(x+a \delta +b\lambda)&=&x+a \delta +b\lambda 
+b\tau +  (\frac{b}{2}(\tau,\tau)-(x,\tau))\delta 
\end{eqnarray}
for $ x \in V$ and $a,b \in \rset$.
Let $S= \{ t_{\tau}; \tau \in Q^{\vee} \}$, then the group $\widehat W$ 
is the semi-direct product of $S$ by $W$.

Consider the $\widehat W$-invariant affine subspace 
$$
E = \{x \in \widehat
V; \ (x,\delta) =1\}= V \oplus \rset\delta + \lambda.
$$ 
Let
$\pi\ :\ E \rightarrow V $ be the  projection $ax+b\delta+\lambda
\mapsto ax$ and  
$$
\begin{array}{lcrcl}
i&:&V& \rightarrow & E \\
&&v & \mapsto & v+\lambda
\end{array}
$$

For $w \in \widehat W$, we set $\overline w = \pi \circ w_{|_E}\circ
i$. The map
$w  \mapsto  \overline w$ defines an injective  morphism of groups from  
$\widehat W$ to $\mathrm{Aut}(V)$. We shall identify $\widehat W$ with
its image $W_{\aff}$ under this
map.

For $\alpha \in \Delta$, 
$\overline {s_{\alpha}}$ is the reflection
$s_{\alpha}$ on $V$ associated to $\alpha$, and for $\tau \in Q^{\vee}$,
$\overline {t_{\tau}}$ is the translation $T_{\tau}$ by the vector $\tau$ on
$V$. For $\alpha \in \Delta^+$, $k \geqslant 0$, $x \in V$, we obtain
that 
$$
\begin{array}{c}
\overline{s_{-\alpha+k\delta}}(x)=x-((x,\alpha)-k)\alpha^{\vee}=T_{k\alpha^{\vee}} \circ
s_{\alpha}(x) \\
\overline{s_{\alpha+k\delta}}(x)=x-((x,\alpha)+k)\alpha^{\vee}=
T_{-k\alpha^{\vee}} \circ s_{\alpha}(x). 
\end{array}
$$
Thus $\overline{s_{-\alpha+k\delta}}$ and
$\overline{s_{\alpha+k\delta}}$ are the orthogonal reflections with
respect to $H_{\alpha,k}=  \{ x
\in V;\ (x,\alpha) =k \}$ and $H_{\alpha,-k}$ respectively. It follows
that $W_{\aff}$ is the semi-direct product
of $W$ by the group of translations $T_{\tau}$, $\tau \in Q^{\vee}$.

Observe that for $v \in W$, $\tau \in Q^{\vee}$, $\alpha \in \Delta$
and $k \in \zset$, we have 
$$
\overline{vt_{\tau}}(H_{\alpha,k})=H_{v(\alpha),k+(\tau,\alpha)}.
$$
Recall that the connected components of the complement in $V$ of
$\bigcup_{\alpha \in \Delta,k \in \zset}H_{\alpha,k}$ are
called alcoves. The group $W_{\aff}$
acts simply transitively on the set of alcoves. We denote 
$$
C=\{x \in V; (\alpha_i,x)
>0 \mbox{ for all } \alpha_i \in \Pi\},\ A=\{x \in C;
(\theta,x)<1\}
$$ 
respectively the fundamental chamber and the
fundamental alcove with respect to $\Pi$ and $\widehat{\Pi}$. 

%---------------------------------------------------------------------------------------------------------------------------------

We shall end this section by recording the following results :

\begin{prop}\label{propdeN}

For $w \in \widehat W$, let $N(w)=\{ \beta \in \widehat{\Delta}^+;\
w^{-1}(\beta)<0\}$ and denote by $\ell (w)$ the length of any reduced
expression of $w$.

{\rm{(a)}} We fix a reduced expression of 
$w = s_{\beta_1} \circ \cdots \circ s_{\beta_k}$ with $ \beta_i
\in \racbase$, then $N(w)= \{s_{\beta_1}\circ\cdots \circ
s_{\beta_{p-1}}(\beta_p);1\leq p\leq  k\}$. 
In particular, 
$N(w)$ contains a simple root.

{\rm{(b)}} Let $w_1,\ w_2 \in \widehat W$, then $N(w_1)
\subseteq N(w_2)$ 
if and only if, there exists $u \in \widehat W$ such that 
$w_2=w_1u$, and $\ell(w_2) = \ell(w_1)
+\ell(u)$. In particular, $w$ is uniquely determined by $N(w)$. 

{\rm{(c)}} If $N(w) \cap \Delta^+\not=\emptyset$, then $N(w)\cap \Pi\not=\emptyset$.

\end{prop}

\begin{proof}
For parts {\rm{(a)}} and {\rm{(b)}}, see for example \cite{CP}. Let us
prove {\rm{(c)}}. The case $\widetilde{A_1}$ is clear. In the others
cases, this is a direct consequence of the fact that $N(w)$ is 
a ``compatible'' set, by theorem 1.3 from \cite{CP}.
\end{proof}

%---------------------------------------------------------------------------------------------------------------
\section{$I$-compatible elements in $\widehat W$}

Let $I \subset \Pi$ and $\ilie$ be an ad-nilpotent ideal of
$\plie_I$.  We set 
$$
\Phi_{\ilie} =  \{ \alpha \in \Delta^+ \setminus \Delta_I;\ 
\mathfrak{g}_{\alpha } \subseteq \ilie \}.
$$ 
Then $\ilie = 
\bigoplus_{\alpha \in \Phi_{\ilie} } \mathfrak{g}_{\alpha}$ and if 
$\alpha \in \Phi_{\ilie}$, $\beta \in \Delta^+\cup \Delta_I$ are such that $\alpha +\beta
\in \Delta^+$, then $\alpha +\beta \in \Phi_{\ilie}$. 

Conversely, set 
$$ 
\calf_I = \{ \Phi \subset \Delta^+\setminus \Delta_I;\mbox{if } \alpha \in \Phi,
\beta \in \Delta^+\cup \Delta_I, \alpha +\beta \in \Delta^+, 
\mbox{then} \ \alpha +\beta\in \Phi \}.
$$
Then for $\Phi \in \calf_I$, $\ilie_{\Phi} = \bigoplus_{\alpha \in
\Phi} \mathfrak{g}_{\alpha}$ is an ad-nilpotent ideal of $\plie_I$. 

We obtain therefore a bijection
$$
\{\mbox{ad-nilpotent ideals of } \plie_I \}  \rightarrow  \calf_I ,\ 
\ilie  \mapsto  \Phi_{\ilie}.
$$
For $\Phi \in \calf_I$, we define $\Phi^1=\Phi$, $\Phi^k=(\Phi^{k-1}+
\Phi) \cap \Delta$, for $k \geq 2$ and 
$$
L_{\Phi}= \bigcup_{k \in \nset^*} (-\Phi^k+k\delta).
$$

Since any ad-nilpotent ideal of $\plie_I$ is an ad-nilpotent ideal of
the Borel subalgebra $\plie_{\emptyset}=\blie$, we have by \cite{CP}
the following proposition :

\begin{prop}

Let $\Phi \in \calf_I$, then there exists a unique $w_{\Phi} \in
\widehat W$ such that $L_{\Phi}=N(w_{\Phi})$.

\end{prop}

Thus we have the following injective map :
$$
\begin{array}{rcl}
\{\mbox{ad-nilpotent ideals of }\plie_I \} & \rightarrow & \widehat W \\
\ilie & \mapsto & w_{\Phi_{\ilie}}
\end{array}
$$

Recall from \cite{CP} the following characterization of the image of
the above map when $I=\emptyset$.

\begin{prop}\label{compatibilité}
Let $w \in \widehat W$, then there exists an ideal $\ilie$ of $\blie$
such that 
$N(w)=L_{\Phi_{\ilie}}$ 
if and only if 

{\rm(a)} $w^{-1}(\alpha) >0$, for all $\alpha \in \Pi$.

{\rm(b)} If $w(\alpha)<0$ for some $\alpha \in \widehat{ \Pi}$, then 
$w(\alpha)=\beta -\delta$ 
for some $\beta \in \Delta^+$.
If these conditions are verified, we say that $w$ is {\bf Borel-compatible}
or {\bf $\emptyset$-compatible}.
\end{prop}

For $w \in \widehat W$, let $\Phi_w=\{ \alpha \in \Delta;
-\alpha+\delta \in N(w)\}$. It follows that if $w$ is
$\emptyset$-compatible, then $\Phi_w \subset \Delta^+$.

\begin{theoreme}\label{theoreme1}
Let $w \in \widehat W$ be Borel-compatible and $I \subset \Pi$. 
The following conditions are equivalent : 

{\rm(a)} $\ilie_{\Phi_w}$ is an ad-nilpotent ideal of $\plie_I$.

{\rm(b)} $s_{\alpha}(\Phi_w)=\Phi_w$, for all $\alpha \in I$.

{\rm(c)} $s_{\alpha}(L_{\Phi_w})=L_{\Phi_w}$, for all $\alpha \in I$.

{\rm(d)} $N(s_{\alpha}w)=N(w) \cup\{\alpha\}$, for all $\alpha \in I$.

{\rm(e)} $w^{-1}(\alpha) \in \widehat{\Pi}$, for all $\alpha \in I$.

\noindent If the hypothesis and these conditions are verified, we say that  $w$
is {\bf $I$-compatible}.

\end{theoreme}

\begin{proof}
$(a)\Rightarrow (b)$
By assumption, we have $\Phi_w \in \calf_I$. Let $\beta \in \Phi_w$,
then $s_{\alpha}(\beta)= 
\beta-(\beta,\alpha^{\vee})\alpha$, hence
$s_{\alpha}(\Phi_w)\subset\Phi_w$, for all $\alpha \in I$. Moreover, 
since $s_{\alpha}$ is an involution, we obtain that 
$s_{\alpha}(\Phi_w)=\Phi_w$.

$(b)\Rightarrow (c)$
Since $w(\delta)=\delta$, for all $w \in
\widehat W$, this is clear (by induction on $k$ or just remark that
${\Phi_w^k} \in \calf_I$).

$(c)\Rightarrow (d)$
Let $\alpha \in I$, by assumption, we have $s_{\alpha}(N(w))=N(w)$, hence for
$\beta \in N(w)$, we have $s_{\alpha}(\beta) \in N(w)$. So
$w^{-1}s_{\alpha}(\beta)<0$ and $\beta \in N(s_{\alpha}w)$. We have 
proved that $N(w) \subset N(s_{\alpha}w)$. Since $\sharp N(w)=\ell (w)$ and
$\ell (s_{\alpha}w)=\ell (w)\pm 1$, by Proposition \ref{propdeN}, we obtain
that $\sharp N(s_{\alpha}w)=\sharp N(w)+1$. Moreover we have 
$(s_{\alpha}w)^{-1}(\alpha)=w^{-1}(-\alpha)<0$, hence $N(s_{\alpha}w)=N(w)
\cup\{\alpha\}$.

$(d)\Rightarrow (e)$
Let $\alpha \in I$. By assumption, we have $N(w) \subset
N(s_{\alpha}w)$, hence by Proposition \ref{propdeN}, there exists
$\beta \in \widehat{\Pi}$ such that, 
$$
N(s_{\alpha}w)=N(ws_{\beta})=N(w)\cup \{w(\beta)\}=N(w) \cup \{\alpha\}.
$$
Consequently, we have $w^{-1}(\alpha)=\beta \in \widehat{\Pi}$.

$(e)\Rightarrow (a)$
Let $\alpha \in I$ and assume that 
$w^{-1}(\alpha) \in \widehat{\Pi}$. Let  $\beta \in \Phi_w$ be such that $\beta- \alpha
\in \Delta^+$. We have 
$$
w^{-1}(-(\beta-\alpha)+\delta)=w^{-1}(-\beta +\delta) +
w^{-1}(\alpha) \in (\widehat{\Delta}^- + \widehat{\Pi})\cap \widehat{\Delta}.
$$
It follows that $w^{-1}(-(\beta-\alpha)+\delta) <0$. Moreover, 
$w^{-1}(-\alpha+\delta)=w^{-1}(-\alpha) +\delta >0$ hence $\alpha
\not \in \Phi_{\ilie_w}$. We obtain that $\Phi_w \in \calf_{\alpha_i}$, for all $\alpha_i
\in I$, hence $\Phi_w$ belongs to $\calf_I$ and $\ilie_{\Phi_w}$ 
is an ideal of $\plie_I$.
\end{proof}

Another characterization of ad-nilpotent ideals in $\blie$ is given in
\cite{CP2} via the set $D=\{\tau \in Q^{\vee}; (\tau, \alpha_j) \leq
1, j=1,\dots, l \mbox{ and }(\tau, \theta) \geq -2\}$. 
Let $\widetilde {D}=\{(\tau,v) \in 
D\times W;vt_{\tau}(A) \subset C\}$. We can state this
characterization in the following way:

\begin{prop}\label{rappelborelsurD}
The following map is bijective : 
$$
\begin{array}{ccl}
\widetilde{D} &\rightarrow &\{w \in \widehat W,\emptyset\mbox{-compatible} \} \\
(\tau,v)& \mapsto &vt_{\tau} .
\end{array}
$$
\end{prop}

\begin{remarque}
In \cite{CP2}, the above correspondence is not viewed in the same way
since the elements of $\widehat W$ are written $t_{\tau}v=vt_{v^{-1}(\tau)}$ 
instead of $vt_{\tau}$, for $w \in W$ and $\tau \in Q^{\vee}$.
\end{remarque}

Let $w \in \widehat W$ be Borel-compatible, then  $I_w=\{ \alpha \in
\Pi; w^{-1}(\alpha) \in \widehat{\Pi}\}$ is the unique maximal element of 
$\{I \subset \Pi;\ w
\mbox{ is } I\mbox{-compatible}\}$. For $\tau \in Q^{\vee}$, set 
$$
D_{\tau}=
\left\{
\begin{array}{ll}
\{\alpha \in
\Pi; (\alpha,\tau)=0 \}\cup\{-\theta\}  &\mbox{ if }(\theta,\tau)=-1, \\
\{\alpha \in\Pi; (\alpha,\tau)=0 \}  &\mbox{ if }(\theta,\tau)\not=-1.
\end{array}
\right.
$$

\begin{prop}\label{condsurtau}
Let $(\tau,v) \in \widetilde D$, and $w=vt_{\tau} \in \widehat
W$. Then $v(D_{\tau})=I_w$. In particular, $w$ is $I$-compatible 
if and only if $I \subset v(D_{\tau})$.
\end{prop}

\begin{proof}
Let $\alpha \in I_w$, then 
$$
w^{-1}(\alpha)=t_{-\tau}v^{-1}(\alpha)=v^{-1}(\alpha)
+(v^{-1}(\alpha),\tau)\delta \in \widehat\Pi.
$$
If $w^{-1}(\alpha) \in \Pi$, then we have $v^{-1}(\alpha) \in \Pi$ and 
$(v^{-1}(\alpha),\tau)=0$, hence $v^{-1}(\alpha)\in D_{\tau}$.  
If $w^{-1}(\alpha)=\alpha_0 $, then we have $v^{-1}(\alpha) =-\theta$ and 
$(\theta,\tau)=-1$, hence $-\theta=v^{-1}(\alpha) \in D_{\tau}$. 

Conversely, let $\alpha \in  D_{\tau} \cap \Pi$, then 
$vt_{\tau}(\alpha)=v(\alpha) \in
\Delta^+$, because $w$ is Borel-compatible. Then we have
$N(ws_{\alpha})=N(w) \cup \{w(\alpha)\}$, and by part {\rm(3)} of 
proposition \ref{propdeN} ,
there exists a simple root $\beta\in \Pi$ such that $\beta \in
N(ws_{\alpha})$. Since $N(w) \cap \Delta^+=\emptyset$, we obtain that
$w(\alpha)=\beta$ and $v(\alpha) \in I_w$.

Assume now that $-\theta\in D_{\tau}$. Since $w$ is Borel-compatible, 
$vt_{\tau}(\alpha_0)=-v(\theta)\in\Delta^+$. As above we have 
$N(ws_{\alpha_0})=N(w) \cup \{w(\alpha_0)\}$, and by part {\rm(3)} of 
proposition \ref{propdeN}, there exists a simple root $\beta\in \Pi$ 
such that $\beta \in N(ws_{\alpha_0})$. Since $N(w) \cap
\Delta^+=\emptyset$, we obtain that
$w(\alpha_0)=\beta$ and $v(-\theta) \in I_w$.

We have therefore proved that $v(D_{\tau})=I_w$, which concludes the proof.
\end{proof}

Let us denote $H_{\alpha}=H_{\alpha,0}$ for $\alpha \in \Pi$, and 
$H_{\alpha_0}=H_{\theta,1}$. Let $\{\omega_1,
  \dots, \omega_l\}$ 
be elements of $V$ such that $(\omega_i,\alpha_j)=\delta_{ij}$. Set
$n_0=1$ and let $n_i$, 
$i=1,\dots, l$, be the strictly positive integers such that 
$\theta =\sum_{i=1}^l n_i\alpha_i$. Let 
$\overline{\omega}_i=\omega_i/n_i$, $i=1\dots l$, and
$\overline{\omega}_0=0$. Then the closure $\overline A$ of $A$ is the
convex hull $\mathrm{Conv}(\overline{\omega}_0,\overline{\omega}_1, \dots, 
\overline{\omega}_l)$ of $\obarre_0, \dots,\obarre_l$. 
For $k \in \nset^*$, the convex hull (resp. the image by 
$\overline w \in W_{\aff}$ of the convex hull) of $(k+1)$ points in $
 \{\overline{\omega}_0,\overline{\omega}_1, \dots, 
\overline{\omega}_l\}$ is called a $k$-face of $\overline A$ (resp. of
$\overline w (\overline A)$). For example, $H_{\alpha_i} \cap
\overline A=\mathrm{Conv}(\obarre_0, \dots,
\obarre_{i-1},\obarre_{i+1}, \dots, \obarre_l)$ is an $(l-1)$-face 
of $\overline A$.

We shall  give yet another characterization of ad-nilpotent ideals of
$\plie_I$ which shall be useful in enumerating abelian ideals when
$\glie$ is of type $A$ or $C$.

\begin{prop}\label{face}
Let $w \in \widehat W$ be Borel-compatible and $I \subset \Pi$. Then, $\ilie_{\Phi_w}$ is
an ideal of $\plie_I$ if and only if for all $\alpha \in I$,
$\overline w(\overline A)\cap H_{\alpha}$ is an $(l-1)$-face of 
$\overline w(\overline A)$.

\end{prop}

\begin{proof}

Assume that $w \in \widehat W$ is $I$-compatible. Let 
$(\tau,v) \in \widetilde D$ be such that $w=vt_{\tau}$. By
proposition \ref{condsurtau}, 
$I\subset v(D_{\tau})$, 
and $v^{-1}(\alpha) 
\in D_{\tau}$, for all $\alpha \in I$. Let $\alpha \in I$, we distinguish two cases :

If $v^{-1}(\alpha)=\beta \in \Pi$, then $(\beta, \tau)=0$. We obtain that 
$$
\overline{vt_{\tau}}(H_{\beta})=H_{v(\beta),(\tau ,\beta )}=H_{\alpha}.
$$
Hence $\overline w(\overline A)\cap H_{\alpha}$ is an $(l-1)$-face of
$\overline w(\overline A)$.

If $v^{-1}(\alpha)=-\theta $, then $(\theta ,
\tau)=-1$. 
We obtain that 
$$
\overline{vt_{\tau}}(H_{\alpha_0})=H_{v(\theta ),(\tau ,\theta  )+1}=H_{\alpha}.
$$
Hence, $\overline w(\overline A)\cap H_{\alpha}$ is an $(l-1)$-face
of $\overline w(\overline A)$.

Conversely, let $v \in W$, $\tau \in Q^{\vee}$ be such that $w=vt_{\tau} \in \widehat W$ 
is Borel-compatible. By assumption, for all $\alpha \in I$, there exists $\beta \in
\widehat{\Pi}$ such that $\overline w(H_{\beta})=H_{\alpha}.$

If $\beta \in \Pi$, then
$$\
\overline{vt_{\tau}}
(H_{\beta})=H_{v(\beta),(\tau ,\beta )}=H_{\alpha}
$$
hence $(\tau,\beta)=0$, and $w^{-1}(\alpha)=\pm \beta$. Since $w$ is 
Borel-compatible, we
have necessarily $w^{-1}(\alpha)>0$, and so $\alpha \in v(D_{\tau}).$

If $\beta=\alpha_0$, then  
$$\overline{vt_{\tau}}(H_{\alpha_0})=
H_{v(\theta ),(\tau ,\theta  )+1}=H_{\alpha}
$$
hence $(\tau,\theta)=-1$, and $w^{-1}(\alpha)=\pm (\theta -\delta)$. 
Since $w$ is Borel-compatible, we
have necessarily $w^{-1}(\alpha)>0$, and so $\alpha \in v(D_{\tau}).$ We have
proved that  $I\subset v(D_{\tau})$, and by
proposition \ref{condsurtau}, $w$ is $I$-compatible.
\end{proof}

Let $H_{\emptyset}=V$. For $J \subset \widehat{\Pi}$ non empty, denote $H_J=\bigcap_{\alpha \in J}
H_{\alpha}$. By the proposition above, if $w$ is $I$-compatible, then
we have $\overline w(\overline A) \cap H_I=\overline
w(\overline A \cap
H_{w^{-1}(I)})$.

%----------------------------------------------------------------------------------------------------------------------------------

\section{Volume of the faces of the fundamental alcove}

Recall from \cite{CP} and \cite{Ko}, that $w \in \widehat W$ is
Borel-compatible and the ideal $\ilie_{\Phi_w}$ of $\blie$ 
is abelian if and only if $\overline w(A) \subset 2A$. As a
consequence, we have the following remarkable result of Peterson : the
number of abelian ideals of $\blie$ is $2^l$. Observe that the above
result says that the number of abelian ideals in $\blie$ depends only
on the rank of $\glie$. In the case of parabolic algebras, we shall
see in this section to what extent this result can be extended. 
 
For $J \subset\widehat{\Pi}$, let 
$F_J=\overline A \cap H_J={\rm{Conv}}(\obarre_j; \ \alpha_j \not\in
J)$. Observe that the $F_J$ are the faces of $\overline A$. 
Let $w \in \widehat W$, if $\overline w(\overline A) \cap H_J$ is
an $(l-\sharp J)$-face of $\overline w(\overline A)$, then we shall
call $\overline w(\overline A) \cap H_J$ an $(l-\sharp
J)$-{\textbf{alcove}} of $H_J$.

\begin{prop}\label{alcove}
{\rm(a)} Let $w \in \widehat W$ and $I \subset \Pi$, if $\overline w (A) \subset 2A$ and
$\overline w(\overline A) \cap H_I$ is an $(l-\sharp I)$-alcove of $H_I$, then 
$w$ is $I$-compatible.

{\rm(b)} Let $I \subset \Pi$ and $w,w' \in \widehat W$ be
$I$-compatible. If $\overline w (A) \subset 2A$, $\overline {w'} (A) \subset 2A$ and
$\overline w(\overline A) \cap H_I=\overline {w'}(\overline A) \cap
H_{I}$, then $w=w'$.
\end{prop}

\begin{proof}
{\rm(a)} Let $w \in \widehat W$ and $I \subset \Pi$ be of cardinality $r$. If $\overline w
(A) \subset 2A$, then $w$ is Borel-compatible and the ideal 
$\ilie_{\Phi_w}$ is abelian. 

Set $N=l-r +1$. Since $\overline w(\overline A) \cap H_I$ is an $(l-r)$-alcove of $H_I$, there exist $N$ vertices $\obarre_{i_1},\dots,\obarre_{i_N}$ of $\overline A$ such that $\overline w(\obarre_{i_1}),\dots, \overline w(\obarre_{i_N})$ belong to $\overline w(\overline A) \cap H_I$.

There exist $r$ distinct reflecting affine hyperplanes $H'_1,\dots, H'_r$ of the form $H_{\alpha}$, for $\alpha\in \widehat{\Pi}$, such that $\bigcap_{j=1}^r H'_j$ contains $\obarre_{i_1},\dots,\obarre_{i_N}$. For $j=1, \dots, r$, $H_I\cap \overline w(H'_j)$ contains $\overline w(\obarre_{i_1}),\dots, \overline w(\obarre_{i_N})$. Since the dimension of $H_I$ is $N-1$, it follows that $H_I\subset\overline w(H'_j)$.

The hyperplane $H_I$ is defined by the equations $(x,\alpha)=0$ for all $\alpha\in I$, it follows that $\overline w(H'_j)$ is an hyperplane of the form $H_{\beta,0}$, where $\beta$ is a linear combination of elements of $I$.

Assume that $\beta\not\in I$. Then, the intersection of $H_{\beta,0}$ with the closure of the fundamental chamber $C$ is of dimension at most $l-2$. Since by construction $H_{\beta,0}$ contains an $(l-1)$-face of $\overline w(\overline A)$, and $\overline w(\overline A)\subset \overline C$, we obtain a contradiction. It follows that $\beta\in I$.

Set $w=vt_{\tau}$. We then have that for each $\beta\in I$ :
$$
w^{-1}(H_{\beta,0})=H_{v^{-1}(\beta),(\tau, v^{-1}(\beta))}=H_{\alpha}
$$
for some $\alpha\in\widehat{\Pi}$. If $\alpha\in\Pi$, then $v^{-1}(\beta)=\pm \alpha$ and $(\tau, v^{-1}(\beta))=0$. Since $w$ is Borel-compatible, we obtain that $w^{-1}(\beta)=\alpha$.

If $\alpha=\alpha_0$, we obtain that $v^{-1}(\beta)=\pm \theta$ and $(\tau, v^{-1}(\beta))=\pm 1$. Since $w$ is Borel-compatible, we finally obtain that $w^{-1}(\beta)=\alpha_0$. Thus, $w^{-1}(I)\subset\widehat{\Pi}$, and $w$ is $I$-compatible as required.

{\rm(b)} Let $I \subset \Pi$ and $w,\ w' \in \widehat W$ be
$I$-compatible. Let $\alpha \in I$, then $w$ is
$I\setminus\{\alpha\}$-compatible. It follows by proposition
\ref{face} that $\overline w (\overline A)
\cap H_{I\setminus\{\alpha\}}$ is an $(l-\sharp I+1)$-alcove of
$H_{I\setminus\{\alpha\}}$ and it is the convex hull of $\overline{w}(\overline A)\cap
H_{I}$ and a vertex of
$H_{I\setminus\{\alpha\}}\cap\overline w(\overline A)$, which is not in
$H_I \cap \overline w(\overline A)$. In the same way $\overline {w'} (\overline A)
\cap H_{I\setminus\{\alpha\}}$ is an $(l-\sharp I+1)$-alcove of
$H_{I\setminus\{\alpha\}}$ and it is the convex hull of $\overline{w'}(\overline A)\cap
H_{I}$ and a vertex of
$H_{I\setminus\{\alpha\}}\cap\overline {w'}(\overline A)$, which is not in
$H_I \cap \overline w(\overline A)$. Since $\overline w(\overline A)
\subset 2\overline A$, there is a unique vertex in
$H_{I\setminus\{\alpha\}}$ satisfying these conditions. So, $\overline w (\overline A)
\cap H_{I\setminus\{\alpha\}}=\overline {w'} (\overline A)
\cap H_{I\setminus\{\alpha\}}$ and by induction, we have 
$\overline{w}(\overline A)=\overline{w'}(\overline
A)$. Hence $w=w'$.
\end{proof}

Let $F'_J=\overline{2A}\cap H_J={\rm{Conv}}(2\obarre_j;\ \alpha_j \not\in
J)$. It is clear that $F'_J$ is a union of $(l-\sharp
  J)$-alcoves of $H_J$. Let 
$$
\mathcal{A}b_I=\{w \in \widehat W; \ilie_{\Phi_w}
\mbox{ is an abelian ideal of } \plie_I\}.
$$
By the above proposition and by proposition \ref{face}, we obtain the following result :

\begin{theoreme}\label{bijectionab}
Let $I \subset \Pi$, then the map $w \mapsto \overline w(\overline A)\cap H_I$ is a bijection between 
$\mathcal {A}b_I$ and the set of all the $(l-\sharp I)$-alcoves of $F'_I$.
\end{theoreme}

\begin{remarque}
The above theorem can be viewed as a generalization of Peterson's
result. 
\end{remarque}

In order to determine $\sharp \mathcal Ab_I$, we are reduced
to computing the volume of the $(l-\sharp I)$-alcoves of $F'_I$. Furthermore, 
to compute the volume of the $(l-\sharp I)$-alcoves of $F'_I$, 
it suffices to compute 
the volume of the $(l-\sharp I)$-faces of $\overline A$. 

Let $d(x,H_{\alpha})$ denote the distance from $x \in V$ 
to the affine hyperplane $H_{\alpha}$, for $\alpha \in
\widehat{\Pi}$. 
For $B$ a $k$-alcove, let
${\rm{Vol}}_k(B)$ be the $k$-volume of $B$. By \cite{Berger}, 
the volume of the fundamental alcove is  
$$
{\rm{Vol}}_l(A)=\frac{1}{l}\times d(0, H_{\alpha_0})\times 
{\rm{Vol}}_{l-1}(F_{\alpha_0}).
$$
Since the projection of $0$ on $H_{\alpha_0}$ is
$\frac{\theta}{\|\theta\|^2}$, we have 
$d(0, H_{\alpha_0})=\frac{1}{\|\theta\|}$.
We obtain that ${\rm{Vol}}_l(A)=\frac{1}{l\|\theta\|}
{\rm{Vol}}_{l-1}(F_{\alpha_0})$. 
Moreover, by \cite{Cox}, 
$$
{\rm{Vol}}_l(A)=\frac{1}{l!}\left|\obarre_1\wedge \cdots
\wedge\obarre_l\right|.
$$
Let $D=\left|\obarre_1\wedge \cdots
\wedge\obarre_l\right|$, 
then 
\begin{equation}\label{volume0}
{\rm{Vol}}_{l-1}(F_{\alpha_0})=\frac{D}{(l-1)!}n_0\|\theta\|.
\end{equation}

To compute the $(l-1)$-volume of the faces $F_{\alpha_i}$, $i=1,
\cdots, l$, we compute the 
$l$-volume of the convex hull of
$\left(\{\obarre_1,\dots,\obarre_l\}\setminus \{\obarre_i\}\right) \cup 
\{\frac{\alpha_i}{\|\alpha_i\|}\}$. 
Thus, we have : 
$$
{\rm{Vol}}_{l-1}(F_{\alpha_i})=\frac{1}{(l-1)!}\left|\obarre_1\wedge \cdots 
\wedge\frac{\alpha_i}{\|\alpha_i\|}\wedge \cdots \wedge\obarre_l\right|.
$$
Since $\alpha_i=\sum_ {k=1}^l (\alpha_i, \alpha_k)\omega_k$, 
\begin{equation}\label{volumei}
{\rm{Vol}}_{l-1}(F_{\alpha_i})=\frac{D}{(l-1)!}n_i\|\alpha_i\|.
\end{equation}

We have therefore computed the $(l-1)$-volume of the $(l-1)$-faces of
$\overline A$. 
In particular, we have :
\begin{lemme}\label{lemmevolume1}
Let $\alpha_i$, $\alpha_j \in \widehat{\Pi}$, be such that $(\alpha_i,\alpha_i)=(\alpha_j,\alpha_j)$, then :
$$
n_i{\rm{Vol}}_{l-1}(F_j)=n_j{\rm{Vol}}_{l-1}(F_i).
$$
\end{lemme}

This Lemma also appears as Proposition $26$ in \cite{Sut}. We shall generalize this result. For $I \subset \widehat{\Pi}$, 
let $n_I=1$ if $I=\emptyset$, and $n_I=\prod_{\alpha_i \in I}n_i$ otherwise. 
We shall prove the following result :

\begin{prop}\label{propvolume}
Let $I \subset \Pi$ and $w \in \widehat W$ be such that $w^{-1}(I)=J
\subset \widehat{\Pi}$. Then, we have :
\begin{eqnarray*}
n_I{\rm{Vol}}_{l-\sharp J}(F_{J})&=&n_J{\rm{Vol}}_{l-\sharp I}(F_{I}) 
\end{eqnarray*}

\end{prop}

%------------------------------------------------------------------------------------------------------

To prove this proposition, we need the following technical lemma :

\begin{lemme}\label{lemmetechnique}
Let $I \subset \Pi$ be such that $\sharp I \leqslant l-1$. Let $w \in
\widehat W$ be such that $w^{-1}(I)=J\subset
\widehat{\Pi}$. Let $\alpha_j$ be any element of $J$ if
$\alpha_0\not\in J$, and $\alpha_j=\alpha_0$ if $\alpha_0\in J$. Set
$\alpha_i=w(\alpha_j)$. Then  we have : 
$$
n_id(\obarre_i,H_I)=n_jd(\obarre_j,H_J).
$$
\end{lemme}

\begin{proof}
The result is clear if $J=\emptyset$. We may therefore 
assume that $1\leqslant \sharp J \leqslant l-1$. 

\medskip
\underline{Step 1} : Assume that $\alpha_0 \in J$. We
shall determine the distance $d(\obarre_0,H_J)$.

Let $J_0$ be the connected component of $J$ containing $\alpha_0$. Set $r=\sharp
J_0$. 

If $J_0=\{\alpha_0\}$, then the projection of $0$ on $H_J$ is
$\frac{\theta}{\|\theta\|^2}$. Therefore, the distance $d(\obarre_0,H_J)$ is
$\frac{1}{\|\theta\|}$. 

Now assume that $J_0\not=\{\alpha_0\}$. Then, $J_0\setminus\{\alpha_0\}$ contains  one or two roots $\beta$ such that
$(\beta,\theta)\not=0$. Set
$J_0=\{\beta_1,\dots,\beta_r\}$, $\alpha_0=\beta_k$ and 
$V_{J_0}=\bigoplus_{\beta_i\in J_0\setminus\{\alpha_0\}} \rset\beta_i$. 

First of all, assume that 
$J_0\setminus\{\alpha_0\}$ contains only one root $\beta_t$ such that
$(\beta_t,\theta)\not=0$. Let $\gamma_t\in V_{J_0}$ be such that
$(\gamma_t,\beta_t)=1$ and $(\gamma_t,\beta_i)=0$ for all $\beta_i\in 
J_0\setminus\{\beta_t,\beta_k\}$. Let $\mu_t=
(\|\theta\|^2(1-\frac{(\gamma_t,\theta)}{2}))^{-1}$ and 
$\beta=\mu_t(\theta-(\theta,\beta_t)\gamma_t)$. Then, we have 
$(\beta,\alpha)=0$ for all  $\alpha \in J_0\setminus \{\alpha_0\}$ and 
$$
(\beta,\theta)=\mu_t[\|\theta\|^2-(\theta,\beta_t)(\gamma_t,\theta)]
=\mu_t\|\theta\|^2[1-\frac{(\gamma_t,\theta)}{2}]
=1.
$$
For all $ x\in H_J$, we have $(\gamma_t,x)=0$, and so 
$$
\begin{array}{ll}
(\beta-x,\beta)&=\mu_t(\theta-(\theta,\beta_t)\gamma_t,\beta-x) \\
&=\mu_t[(\theta,\beta)-(\theta,x)] \\
&=0.
\end{array}
$$

We have proved that $\beta$ is the projection of $\obarre_0$ in
$H_J$. It follows that by taking any $x\in H_J$, we have $d(\obarre_0,H_J)^2=\|\beta\|^2=(x,\beta)=\mu_t(x,\theta)=\mu_t.$

Since $I\subset \Pi$ and $J=w^{-1}(I)$ and $I$ have the same Dynkin diagram, we have by a case by case consideration that $J_0$ is of type $A_r$, $C_r$, or $D_r$. 

If $J_0$ is of type $A_r$, then by renumbering the roots $\beta_i$, 
the Dynkin diagram of $J_0$ is of the
form :
$$
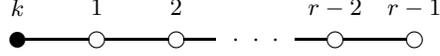
\begin{Dynkin}
\Dbloc{\Dbullet\Deast\Dtext{t}{k}}
\Dbloc{\Dcirc\Dwest\Deast\Dtext{t}{1}}
\Dbloc{\Dcirc\Dwest\Deast\Dtext{t}{2}}
\Dbloc{\Ddots}
\Dbloc{\Dcirc\Dwest\Deast\Dtext{t}{r-2}}
\Dbloc{\Dcirc\Dwest\Dtext{t}{r-1}}
\end{Dynkin}
$$
Then $t=1$, and take 
$$
\gamma_t=\frac{2}{r\|\beta_1\|^2}((r-1)\beta_1+(r-2)\beta_2+\dots
+\beta_{r-1}).
$$
So $(\gamma_t,\theta)=\frac{r-1}{r}$, and we have
\begin{equation}\label{distA_r,0}
\mu_t=\frac{2r}{(r+1)\|\theta\|^2}.
\end{equation}

If $J_0$ is of type $C_r$, then the Dynkin diagram of $J_0$ is of the
form :
$$
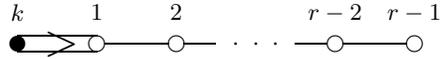
\begin{Dynkin}
\Dbloc{\Dbullet\Ddoubleeast\Dtext{t}{k}}
\Drightarrow
\Dbloc{\Dcirc\Ddoublewest\Deast\Dtext{t}{1}}
\Dbloc{\Dcirc\Dwest\Deast\Dtext{t}{2}}
\Dbloc{\Ddots}
\Dbloc{\Dcirc\Dwest\Deast\Dtext{t}{r-2}}
\Dbloc{\Dcirc\Dwest\Dtext{t}{r-1}}
\end{Dynkin}
$$
Again $t=1$, and take 
$$
\gamma_t=\frac{2}{r\|\beta_1\|^2}((r-1)\beta_1+(r-2)\beta_2+\dots
+\beta_{r-1}).
$$
So $(\gamma_t,\theta)=\frac{2(r-1)}{r}$, and we have
\begin{equation}\label{distC_r,0}
\mu_t=\frac{r}{(r+1)\|\theta\|^2}.
\end{equation}

If $J_0$ is of type $D_r$, then the Dynkin diagram of $J_0$ is of the
form :
$$
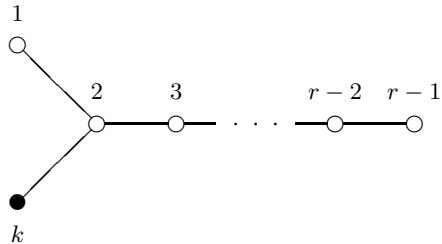
\begin{Dynkin}
\Dbloc{\Dcirc\Dsoutheast\Dtext{t}{1}}
\Dskip
\Dspace\Dbloc{\Dcirc\Dnorthwest\Dsouthwest\Deast\Dtext{t}{2}}
\Dbloc{\Dcirc\Dwest\Deast\Dtext{t}{3}}
\Dbloc{\Ddots}
\Dbloc{\Dcirc\Dwest\Deast\Dtext{t}{r-2}}
\Dbloc{\Dcirc\Dwest\Dtext{t}{r-1}}
\Dskip
\Dbloc{\Dbullet\Dnortheast\Dtext{b}{k}}
\end{Dynkin}
$$
or of the form :
$$
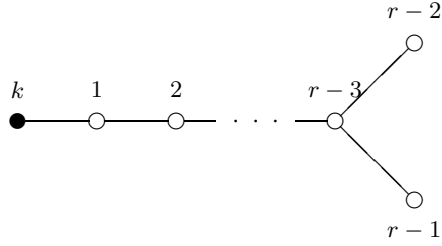
\begin{Dynkin}
\Dbloc{}\Dbloc{}\Dbloc{}\Dbloc{}\Dbloc{}\Dbloc{\Dcirc\Dsouthwest\Dtext{t}{r-2}}
\Dskip
\Dbloc{\Dbullet\Deast\Dtext{t}{k}}
\Dbloc{\Dcirc\Dwest\Deast\Dtext{t}{1}}
\Dbloc{\Dcirc\Dwest\Deast\Dtext{t}{2}}
\Dbloc{\Ddots}
\Dbloc{\Dcirc\Dwest\Dnortheast\Dsoutheast\Dtext{t}{r-3}}
\Dskip
\Dbloc{}\Dbloc{}\Dbloc{}\Dbloc{}\Dbloc{}\Dbloc{\Dcirc\Dnorthwest\Dtext{b}{r-1}}
\end{Dynkin}
$$

In the first case we have $t=2$, and we take 
$$
\gamma_t=\frac{2}{r\|\beta_2\|^2}((r-2)\beta_1+2(r-2)\beta_2+2(r-3)\beta_3+\dots
+2\beta_{r-1}).
$$

Thus $(\gamma_t,\theta)=\frac{2(r-2)}{r}$, and we have
\begin{equation}\label{distD_r,0}
\mu_t=\frac{r}{2\|\theta\|^2}.
\end{equation}

In the second case, we have $t=1$ and we take 
$$
\gamma_t=\frac{1}{\|\beta_1\|^2}(2\beta_1+2\beta_2+\dots+2\beta_{r-3}+\beta_{r-2}+\beta_{r-1}).
$$
Thus we have
\begin{equation}\label{distD_5,0}
\mu_t=\frac{2}{\|\theta\|^2}.
\end{equation}

Assume now that $J_0$ contains two roots $\alpha$ 
such that $(\alpha,\theta)\not=0$. Then the Dynkin diagram of $J_0$ 
is of type $A_r$ and
these two roots are $\beta_{k-1},\beta_{k+1}$ :
$$
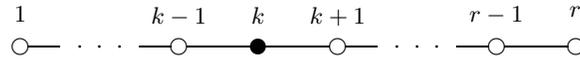
\begin{Dynkin}
\Dbloc{\Dcirc\Deast\Dtext{t}{1}}
\Dbloc{\Ddots}
\Dbloc{\Dcirc\Dwest\Deast\Dtext{t}{k-1}}
\Dbloc{\Dbullet\Dwest\Deast\Dtext{t}{k}}
\Dbloc{\Dcirc\Dwest\Deast\Dtext{t}{k+1}}
\Dbloc{\Ddots}
\Dbloc{\Dcirc\Dwest\Deast\Dtext{t}{r-1}}
\Dbloc{\Dcirc\Dwest\Dtext{t}{r}}
\end{Dynkin}
$$
Let $\eta,\eta'\in V_{J_0}$ be such that
$(\eta,\beta_{k-1})=1=(\eta',\beta_{k+1})$ and
$(\eta,\beta_i)=0$ (resp. $(\eta',\beta_i)=0$) 
for all $\beta_i\in 
J_0\setminus\{\beta_{k-1},\beta_k\}$ (resp. $\beta_i\in 
J_0\setminus\{\beta_k,\beta_{k+1}\}$). Let 
$\mu= (\|\theta\|^2(1-\frac{(\eta+\eta',\theta)}{2}))^{-1}$ and 
$\beta=\mu(\theta-((\theta,\beta_{k-1})\eta+(\theta,\beta_{k-1})\eta'))$. Then we have 
$(\beta,\alpha)=0$  for all  $\alpha \in J_0\setminus\{\alpha_0\}$ and  
$$
\begin{array}{c}
(\beta,\theta)=\mu[\|\theta\|^2
-((\theta,\beta_{k-1})\eta+(\theta,\beta_{k-1})\eta',\theta)]=1 \\
(\beta-x,\beta)=\mu(\theta-
((\theta,\beta_{k-1})\eta+(\theta,\beta_{k-1})\eta')
,\beta-x)=0
\end{array}
$$
for all $ x\in H_J$. We obtain that $\beta$ is the projection of $0$ on $H_J$.
Take
\begin{eqnarray*}
\eta&=&\frac{2}{k\|\beta_{k-1}\|^2}((k-1)\beta_{k-1}+(k-2)\beta_{k-2}+\dots
+\beta_{1})\\
\eta'&=&\frac{2}{(r-k+1)\|\beta_{k+1}\|^2}(\beta_{r}+2\beta_{r-1}+\dots
+(r-k)\beta_{k+1})\\
\end{eqnarray*}
then $(\eta+\eta',\theta)=\frac{k-1}{k}+\frac{r-k}{r-k+1}$. We obtain that :
\begin{equation}\label{dist2A_r,0}
d^2(\obarre_0,H_J)=\|\beta\|^2=\mu=\frac{2k(r-k+1)}{n_0^2(r+1)\|\theta\|^2}.
\end{equation}

Observe that the formulas \eqref{distA_r,0} and \eqref{dist2A_r,0}
generalize the formula obtained when $J_0=\{\alpha_0\}$. 
Let $k$ be the position of $\alpha_0\in J_0$, then we can sum up the
above results in the following table,
when $1\leqslant \sharp J\leqslant l-1$ :

\begin{table}[h]
$$
\begin{array}{|*{6}{c|}}
\hline
\displaystyle J_0 & \displaystyle A_r 
& \begin{array}{c} C_r \end{array} 
& \begin{array}{c} D_r \\ t=2 \end{array} 
& \begin{array}{c} D_r \\ t=1 \end{array}\\
\hline 

\displaystyle d(\obarre_0,H_J)^2 
& \displaystyle\frac{2k(r-k+1)}{n_0^2(r+1)\|\theta\|^2}
& \displaystyle\frac{r}{n_0^2(r+1)\|\theta\|^2}
& \displaystyle\frac{r}{2n_0^2\|\theta\|^2} 
& \displaystyle\frac{2}{n_0^2\|\theta\|^2}\\
\hline
\end{array}
$$
\caption{\label{tableau1}}
\end{table}

\medskip
\underline{Step 2} : Assume that $J \subset \Pi$. Let $\alpha_j \in
J$. We shall determine the distance $d(\obarre_j,H_J)$.

We have 
$H_J ={\rm{Vect}}(\obarre_t;t\mbox{ such that }\alpha_t\not\in J)\subset
H_{J\setminus\{\alpha_j\}} \subset V$. Let $H_J^{\bot}=\{ x\in V;\ 
(x,\obarre_t)=0\mbox{ for all $t$ such that }
\alpha_t\not\in J\}$, then $H_J^{\bot}={\rm{Vect}}(\alpha_t;\ \alpha_t\in J)$, and 
${\rm{dim}}(H_J^{\bot} \cap H_{J\setminus\{\alpha_j\}})=1$. Since 
$$
H_{J\setminus\{\alpha_j\}} \cap  H_J^{\bot}  =\left\{ x=\sum_{\alpha_t\in
  J}\tau_t\alpha_t;\ (x,\beta)=0 \mbox{ for all } \beta 
\in J\setminus\{\alpha_j\}\right\}, 
$$
there exists $\gamma \in V$ such that 
$H_{J\setminus\{\alpha_j\}} \cap  H_J^{\bot}={\rm{Vect}}(\gamma) $, 
and $(\gamma,\alpha_j)\not=0$. 
Thus, we have $\ H_{J\setminus\{\alpha_j\}}=H_J\oplus \cset\gamma$. It
follows that there exists  $\mu \in \cset^*$ such that 
$\obarre_j+\mu\gamma \in H_J$. 

Let $J_j$ be the connected component of $J$ 
which contains  $\alpha_j$. Set $J_j=\{\beta_1,\dots,\beta_r\}$, 
with $\alpha_j=\beta_k$. Set $V_{J_j}=\bigoplus_{\beta_j \in
  J_j}\rset\beta_j$. We may choose $\gamma$ such that $(\gamma,\alpha_j)=1$
and $(\gamma,\alpha)=0$ for all $\alpha\in J_j\setminus\{\alpha_j\}$.
Since $(\obarre_j
+\mu\gamma,\alpha_j)=0$, we obtain that $\mu=-\frac{1}{n_j}$, and 
hence 
\begin{equation}\label{distance}
d(\obarre_j,H_J)=\frac{\|\gamma\|}{n_j},
\end{equation}
where $\gamma$ depends only on the position of $\alpha_j$ in the Dynkin 
diagram of $J_j$. 

Finally, we need to compute explicitely 
$d(\obarre_j,H_J)$ in some particular cases. We use the numbering of 
\cite[chapter 18]{TY}.

If $J_j=A_r$, take
\begin{multline*}
\gamma=\frac{2}{(r+1)\|\beta_k\|^2}[((r-k+1)\beta_1+
2(r-k+1)\beta_2+\cdots \\ +(k-1)(r-k+1)
\beta_{k-1} 
+k(r-k+1)\beta_k+k(r-k)\beta_{k+1}+\cdots+k\beta_r].
\end{multline*}

If $J_j=C_r$, and $k=r$, take 
$$
\gamma=\frac{2}{\|\beta_r\|^2}(\beta_{1}+2\beta_2+\cdots+(r-1)\beta_{r-1}+
\frac{r}{2}\beta_r).
$$

If $J_j=D_r $, take 
$$
\begin{array}{rcll}
\gamma&=&
\displaystyle\frac{1}{\|\beta_r\|^2}[\beta_{1}+2\beta_2+\cdots+(r-2)\beta_{r-2}
  & \\
&&\displaystyle +\frac{1}{2}[(r-2)\beta_{r-1}+r\beta_r]] & \mbox{if } k=r,\\ 

\gamma&=&
\displaystyle\frac{1}{\|\beta_r\|^2}[\beta_{1}+2\beta_2+\cdots+(r-2)\beta_{r-2} &\\
&&\displaystyle +\frac{1}{2}[r\beta_{r-1}+(r-2)\beta_{r}]] &\mbox{if } k=r-1, \\

\gamma&=&
\displaystyle\frac{1}{\|\beta_1\|^2}[2\beta_{1}+2\beta_2+\cdots+2\beta_{r-2} &\\
&&+\beta_{r-1}+\beta_{r}] &\mbox{if } k=1.
\end{array}
$$

In these particular cases, we obtain the following result : 
\begin{table}[h]
$$
\begin{array}{|*{5}{c|}}
\hline
\displaystyle J_j & \displaystyle A_r 
& \begin{array}{c} C_r \\k=r\end{array} 
& \begin{array}{c} D_r \\k=r-1,r\end{array} 
& \begin{array}{c} D_r \\k=1\end{array}\\
\hline 

\displaystyle d(\obarre_j,H_J)^2 
& \displaystyle\frac{2k(r-k+1)}{n_j^2(r+1)\|\alpha_j\|^2}
& \displaystyle\frac{r}{n_j^2(r+1)\|\alpha_j\|^2}
& \displaystyle\frac{r}{2n_j^2\|\alpha_j\|^2} 
& \displaystyle\frac{2}{n_j^2\|\alpha_j\|^2} \\
\hline
\end{array}
$$
\caption{\label{tableau2}}
\end{table}

\underline{Final step} : We are now in a position to prove the lemma. 
Let $I_i$ be the connected component of $I$ containing $\alpha_i$. 
If $J\subset \Pi$, then we
have the result by \eqref{distance}, since
$\alpha_j$ and $\alpha_i$ have the same position in the Dynkin diagram
of $J_j $ and $I_i$ respectively.

If $\alpha_0\in J$, then the connected component $J_0$ of $J$
containing $\alpha_0$ is of the type $A_r$, $C_r$, or $D_r$. Again since $w^{-1}(\alpha_0)$ and $\alpha_0$ have the same position in the respective Dynkin diagram, we obtain the result by inspecting the 
correspondence between
tables \ref{tableau1} and \ref{tableau2}.
\end{proof}

\begin{proof}[Proof of proposition \ref{propvolume}]
The case $\sharp I=0$ is trivial since in this case,
$F_I=F_J=\overline A$. 

Let $I \subset \Pi$. Let us proceed by induction on $\sharp I$. 
If $\sharp I=1$, the result is proved in lemma
\ref{lemmevolume1}. Assume that
$l>\sharp I >1$ 
and that the claim is true for $\sharp I-1$. Let $\alpha_j$ be any element of $J$ if
$\alpha_0\not\in J$, and $\alpha_j=\alpha_0$ if $\alpha_0\in J$. Set
$\alpha_i=w(\alpha_j)$. Then, we have by lemma \ref{lemmetechnique},
\begin{eqnarray*}
n_J {\rm{Vol}}_{l-\sharp I}(F_{I})& = & n_J(l-\sharp I+1){\rm{Vol}}_{l-\sharp
  I+1}(F_{I\setminus\{\alpha_i\}})\times \frac{1}{d(\obarre_i,H_I)} \\
&=&n_j(l-\sharp I+1)\frac{n_I}{n_i}{\rm{Vol}}_{l-\sharp
  I+1}(F_{J\setminus\{\alpha_j\}})\times \frac{n_i}{n_jd(\obarre_j,H_J)} \\
&=&n_I(l-\sharp I+1){\rm{Vol}}_{l-\sharp
  I+1}(F_{J\setminus\{\alpha_j\}})\times \frac{1}{d(\obarre_j,H_J)} \\
&=&n_I {\rm{Vol}}_{l-\sharp I}(F_{J}).
\end{eqnarray*}
Finally, the result is clear if $\sharp I=l$ since in this case
$F_I$ (resp. $F_J$) is a single point.
\end{proof}

Observe that for $I\subset \Pi$, $F'_I=2F_I$, so 
\begin{eqnarray}\label{eqvolume}
{\rm{Vol}}_{l-\sharp I}(F'_{I})=2^{l-\sharp I}{\rm{Vol}}_{l-\sharp I}(F_{I}). 
\end{eqnarray}

%---------------------------------------------------------------------------------------------------------
We obtain a generalization of Peterson's result :

\begin{theoreme}\label{thabelien}
Let $I\subset \Pi$, then 
$$
\frac{1}{n_I}\sum_{w\in\mathcal Ab_I} 
n_{w^{-1}(I)}=2^{l-\sharp I}. 
$$
\end{theoreme}

\begin{proof}
Let $I \subset \Pi$ and $w \in \widehat W$. By propositions \ref{face} 
and \ref{alcove}, then
$$
\sum_{w\in \mathcal Ab_I}
{\rm{Vol}}_{l-\sharp I}(\overline w(\overline A)\cap H_I)={\rm{Vol}}_{l-\sharp I}(F'_I).
$$
Observe that for $w \in \mathcal Ab_I$ we have $\overline w(\overline A)\cap H_I=
\overline w(F_{w^{-1}(I)})$, and ${\rm{Vol}}_{l-\sharp I}(\overline w(F_{w^{-1}(I)}))
={\rm{Vol}}_{l-\sharp I}(F_{w^{-1}(I)})$. So by proposition \ref{propvolume} 
 and by \eqref{eqvolume}, we obtain that 
\begin{equation*}\label{eqabelien}
\sum_{w \in \mathcal Ab_I}
\frac{n_{w^{-1}(I)}}{n_I}{\rm{Vol}}_{l-\sharp I}(F_{I})=2^{l-\sharp I
}{\rm{Vol}}_{l-\sharp I}(F_I)  
\end{equation*}
Thus, we have the result.
\end{proof}

\begin{theoreme}\label{theoremeab}
Let $I \subset \Pi$, if $\glie$ is of type $A_l$ or $C_l$, then 
the parabolic subalgebras  
$\plie_I$ have exactly $2^{l-\sharp I}$ abelian ideals.

\end{theoreme}

\begin{proof}
If $\glie$ is of type $A_l$ or $C_l$, the numbers $n_i$, for $i=0,\dots,l$, depends only on the length
of $\alpha_i$. It follows that for any $w\in\mathcal Ab_I$,
$n_I=n_{w^{-1}(I)}$. So by theorem \ref{thabelien}, we obtain the result. 
\end{proof}

\begin{remarque}
The fact that the integers $n_i$, for $i=0,\dots,l$, depends only on the length of $\alpha_i$ is false
when $\glie$ is not of type $A$ or $C$. Indeed, theorem \ref{theoremeab} is false in
general. For example, in $B_3$, the parabolic 
subalgebra $\plie_{\{\alpha_1\}}$ has only $3$ abelian ideals. We
shall see in the next section another way to count the number of
abelian ideals in cases $B$ and $D$.
\end{remarque}

%------------------------------------------------------------------------------------------------------

\section{Enumeration of ideals via diagrams}

In this section, we shall determine, via diagram enumeration, the number of ad-nilpotent
(resp. abelian) ideals of $\plie_I$, for $I\subset \Pi$, when $\glie$ 
is simple and of classical type. We shall use the numbering of simple
roots of \cite[chap.18]{TY}.

Recall the following partial order on $\Delta^+$ : 
$\alpha \leqslant \beta$ if $\beta -\alpha$ is a 
sum of positive roots. Then it is easy to see that $\Phi \in
\calf_{\emptyset}$ 
if and only if for all $\alpha \in \Phi, \beta \in \Delta^+$, such
that $\alpha \leqslant \beta$, then 
$\beta \in \Phi$. When $\glie $ is of type 
$A$, $B$, $C$ or $D$, we can display the positive 
roots into a diagram of suitable shape, as in \cite{CP}. Then, 
they established 
a bijection between elements of 
$\calf_{\emptyset}$ and certain subdiagrams.

Let $I\subset \Pi$. In order to adapt this
construction in the parabolic case $\plie_I$, we shall use a 
similar construction, 
but our diagram will depend not only on the type of $\glie$, but also on 
$I$. 

Let $I\subset \Pi$ and $\gamma,\beta\in \Delta^+$. We say that $\beta
\stackrel{I}{\rightarrow}\gamma$ if there exists $\eta\in I$ such that
$\beta+\eta=\gamma$. Define an equivalence relation on 
$ \Delta^+\setminus\Delta_I$ : for $I\subset \Pi$, $\gamma \sim_I \beta$ if there exist
$\beta_1,\dots,\beta_s \in \Delta^+\setminus \Delta_I$ such that 
$$
\begin{array}{rl}
{\rm (i)} &\beta=\beta_1,\ \gamma=\beta_s, \\
{\rm (ii)} & \mbox{either } \beta_i
\stackrel{I}{\rightarrow}\beta_{i+1} \mbox{ or } \beta_{i+1}
\stackrel{I}{\rightarrow}\beta_{i},\ \mbox{ for } i=1,\dots,s-1.
\end{array}
$$
As the standard Levi factor of $\plie_I$ acts in a reductive way on
the nilpotent radical, the fact that two roots $\beta,\gamma$ are $\sim_I$
equivalent means that $\glie_{\alpha}$ and $\glie_{\beta}$ are in the
same simple submodule.

Let $X$ be the type of $\glie$. The idea is to start by displaying 
the positive roots $\Delta^+$ in a diagram $T_X$ of a suitable shape as in \cite{CP} : 
that is, we assign to each box labelled $(i,j)$ in $T_X$, 
a positive root $t_{i,j}$. The shape and the filling of $T_X$ are
chosen such that we obtain a bijection 
between elements of $\calf_{\emptyset}$ and the northwest flushed 
subdiagrams, henceforth nw-diagrams, of $T_X$ (in type $D$, we need to
include also nw-diagrams modulo a permutation of certain columns). 
Then, for $I \subset \Pi$, we delete the boxes containing
elements of $\Delta_I$. Observe that the set of boxes of the same equivalent 
class is connected. Therefore, we can regroup into a big box all the roots 
of the same equivalent class. We obtain a new diagram denoted by $T_X^I$. 
Then, we count the nw-diagrams of $T_X^I$ (again in type $D$, we need to
count also nw-diagrams modulo a permutation of certain columns), 
which are clearly in bijection with the elements of $\calf_I$.

\subsection{Type $A_l$}\label{sectionA}

If $\glie$ is of type $A_l$, then $T_{A_l}$ is a diagram 
of shape $[l,l-1,\dots,1]$. The label $(i,j)$ means a box in the $i$-th
row and the $j$-th column. The boxes $(i,j)$ of $T_{A_l}$ are filled by 
the positive roots $t_{i,j}=\alpha_i+\cdots +\alpha_{l-j+1}$, $1\leqslant i,
j\leqslant l$. For example, for $l=5$, we have :
$$
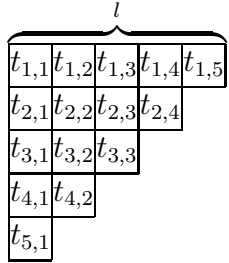
\overbrace{
\vbox{\offinterlineskip
\halign{#&#&#&#&#&#\cr
\htrait&\htrait&\htrait&\htrait&\htrait\cr
\sboite{$t_{1,1}$}&\sboite{$t_{1,2}$}&\sboite{$t_{1,3}$}&\sboite{$t_{1,4}$}&%
\sboite{$t_{1,5}$}&\vtrait\cr
\htrait&\htrait&\htrait&\htrait&\htrait\cr
\sboite{$t_{2,1}$}&\sboite{$t_{2,2}$}&\sboite{$t_{2,3}$}%
&\sboite{$t_{2,4}$}&\vtrait\cr
\htrait&\htrait&\htrait&\htrait\cr
\sboite{$t_{3,1}$}&\sboite{$t_{3,2}$}&\sboite{$t_{3,3}$}&\vtrait\cr
\htrait&\htrait&\htrait\cr
\sboite{$t_{4,1}$}&\sboite{$t_{4,2}$}&\vtrait\cr
\htrait&\htrait\cr
\sboite{$t_{5,1}$}&\vtrait\cr
\htrait\cr
}
}
}^{l}
$$
Let $I\subset \Pi$. We first delete 
the boxes containing elements of $\Delta_I$. 
Then, we regroup the equivalent classes of $\sim_I$ proceeding simple 
root by simple root : for each $\alpha_i\in I$, we regroup the
$(l-i+1)$-th and the $(l-i+2)$-th columns if $i\not=1$, and the rows
$i$, $i+1$ if $i\not=l$, on $T_{A_l}$. 
At the end, we obtain that 
$T_{A_l}^I$ is a  diagram 
of shape 
$[l-\sharp I,l-\sharp I-1, \dots,1]$. 
For example, for $A_5$ and $I=\{\alpha_2, \alpha_3\}$, we have :
$$
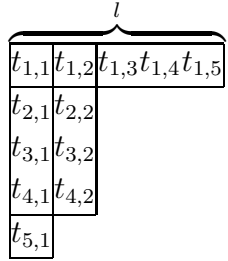
\overbrace{
\vbox{\offinterlineskip
\halign{#&#&#&#&#&#\cr
\htrait&\htrait&\htrait&\htrait&\htrait\cr
\sboite{$t_{1,1}$}&\sboite{$t_{1,2}$}&\sboite{$t_{1,3}$}&\hboite{$t_{1,4}$}&
\hboite{$t_{1,5}$}&\vtrait\cr
\htrait&\htrait&\htrait&\htrait&\htrait\cr
\sboite{$t_{2,1}$}&\sboite{$t_{2,2}$}&\vtrait\cr
\sboite{$t_{3,1}$}&\sboite{$t_{3,2}$}&\vtrait\cr
\sboite{$t_{4,1}$}&\sboite{$t_{4,2}$}&\vtrait\cr
\htrait&\htrait\cr
\sboite{$t_{5,1}$}&\vtrait\cr
\htrait\cr
}
}
}^{l}
$$

Let $\mathcal C_l=\displaystyle \frac{1}{l+1}
\left(\begin{array}{c} 
2l \\l
\end{array}\right)$ denote the $l$-th Catalan number.

\begin{prop}\label{enumerationA}
Let $I\subset\Pi$. 
Let $\mathcal S_{A_l}^I$ be the set of all nw-diagrams of $T_{A_l}^I$. 
Then, the cardinality of $\mathcal S_{A_l}^I$ is $\mathcal C_{l-\sharp I+1}.$
\end{prop}

\begin{proof}
Let $I\subset \Pi$, then $T_{A_l}^I$ if of shape 
$[l-\sharp I,l-\sharp I-1, \dots,1]$, so by \cite[8,6.19 vv.]{Sa}, 
we obtain that the cardinality of the set of nw-diagrams of 
$T_{A_l}^I$ is $\mathcal C_{l-\sharp I+1}.$
\end{proof}

\subsection{Type $C_l$}\label{sectionC}

\begin{definition}
Let $p,q$ be two integers such that $q \leqslant p$. Let $T_{p,q}$ be the
(shifted) diagram of shape $[p+q-1,p+q-3,\dots,p-q+1]$ arranged in the
following way :
$$
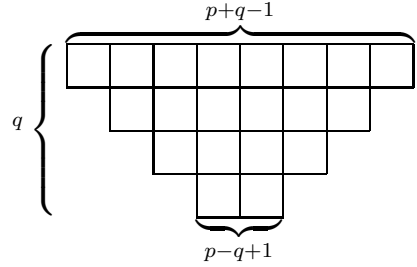
\vbox{\offinterlineskip
\halign{#\cr
${}^{q}\left\{ \vrule height36pt width0pt\right.$\cr
\vrule height1pt width0pt\cr
}
}
\overbrace{
\vbox{\offinterlineskip
\halign{#&#&#&#&#&#&#&#&#\cr
\htrait&\htrait&\htrait&\htrait&
\htrait&\htrait&\htrait&\htrait&\cr
\sboite{}&\sboite{}&\sboite{}&\sboite{}&
\sboite{}&\sboite{}&\sboite{}&\sboite{}&\vtrait\cr
\htrait&\htrait&\htrait&\htrait&
\htrait&\htrait&\htrait&\htrait&\cr
&\sboite{}&\sboite{}&\sboite{}&\sboite{}&
\sboite{}&\sboite{}&\sboite{}&\cr
&\htrait&\htrait&\htrait&\htrait&\htrait&\htrait\cr
&&\sboite{}&\sboite{}&\sboite{}&\sboite{}&\vtrait\cr
&&\htrait&\htrait&\htrait&\htrait\cr
&&&\sboite{}&\sboite{}&\vtrait\cr
&&&\multispan2\hrulefill&\cr
&&&\multispan2$\underbrace{\hboite{}\hboite{}}_{p-q+1}$&\cr
}
}
}^{p+q-1}
$$
\end{definition}

If $\glie$ is of type $C_l$, then $T_{C_l}$ is the diagram $T_{l,l}$, 
and the boxes $(i,j)$ of $T_{C_l}$ are filled by 
the positive roots $t_{i,j}$, where
$$
t_{i,j}=\left\{
\begin{array}{ll}
\alpha_i+\cdots +\alpha_{j-1}+2(\alpha_{j}+\cdots+\alpha_{l-1})+\alpha_l,&
1\leqslant j\leqslant l-1, \\
\alpha_i+\cdots+\alpha_{2l-j},& l\leqslant j\leqslant 2l-1.
\end{array}
\right.
$$
Let $I\subset \Pi$, we first delete 
the boxes containing elements of $\Delta_I$. 
Then, we regroup the equivalent classes of $\sim_I$ proceeding simple 
root by simple root : for each 
$\alpha_i\in I\setminus\{\alpha_l\}$, we first regroup column $2l-i$ 
and column $2l-i+1$ if $i\not=1$, then we regroup the $i$-th and $(i+1)$-th columns
and also the $i$-th and $(i+1)$-th 
rows on $T_{C_l}$. If $\alpha_l\in I$, we regroup also the 
columns $l$ and $l+1$. We obtain at the end that 
$T_{C_l}^I$ is a  diagram 
of shape $T_{l-\sharp I,l-\sharp I}$ if $\alpha_l\not\in I$ and of
shape $T_{l-\sharp I+1,l-\sharp I}$, if $\alpha_l \in I$. 

By \cite{Pr}, we obtain directly that the number of nw-diagram 
of $T_{p,q}$ is $\left(\begin{array}{c}p+q \\ p\end{array}\right)$. 
Consequently, we have the following proposition :

\begin{prop}\label{enumerationC}
Let $I\subset\Pi$. 
Let $\mathcal S_{C_l}^I$ be the set of all nw-diagrams of
$T_{C_l}^I$. Then, the cardinality of $\mathcal S_{C_l}^I$ is
$$
(l-\sharp I +1)\mathcal C_{l-\sharp I} \ \mbox{ if }\  \alpha_l\not\in I,
\quad\mbox{and}\quad
\frac{l-\sharp I +2}{2}\mathcal C_{l-\sharp I+1} 
\ \mbox{ if }\  \alpha_l\in I.
$$
\end{prop}

%------------------------------

\subsection{Type $B_l$ and $D_l$}

Let $I\subset \Pi$. Assume that $\glie$ is of type 
$X=B_l$ or $D_l$. Then the shape of $T_X^I$ is more complicated 
than in the case $A$ or $C$, so we need more combinatorial results 
on diagrams.

\begin{definition}
Let $p,q$ be two integers such that $q \leqslant p$. Let $T'_{p,q}$ be the 
diagram of $q$ rows of the shape $[p,p-1,\dots,p-q+1]$ arranged in the 
following way:
$$
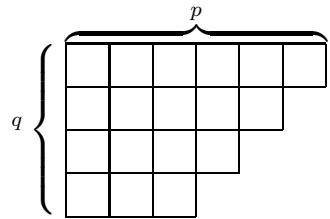
\vbox{\offinterlineskip
\halign{#\cr
${}^{q}\left\{ \vrule height36pt width0pt\right.$\cr
\vrule height1pt width0pt\cr
}
}
\overbrace{
\vbox{\offinterlineskip
\halign{#&#&#&#&#&#&#\cr
\htrait&\htrait&\htrait&\htrait&
\htrait&\htrait\cr
\sboite{}&\sboite{}&\sboite{}&\sboite{}&
\sboite{}&\sboite{}&\vtrait\cr
\htrait&\htrait&\htrait&\htrait&
\htrait&\htrait\cr
\sboite{}&\sboite{}&\sboite{}&\sboite{}&
\sboite{}&\vtrait\cr
\htrait&\htrait&\htrait&\htrait&\htrait\cr
\sboite{}&\sboite{}&\sboite{}&\sboite{}&\vtrait\cr
\htrait&\htrait&\htrait&\htrait\cr
\sboite{}&\sboite{}&\sboite{}&\vtrait\cr
\htrait&\htrait&\htrait\cr
}
}
}^{p}
$$
\end{definition}

\begin{prop}\label{compteA}
Let $p,q$ be two integers such that $q \leqslant p$. Then, 
the number of nw-diagrams of $T'_{p,q}$ is 
$$
\mathcal T'_{p,q}=\frac{(p+q+1)!(p-q+2)}{q!(p+2)!}.
$$
\end{prop}

\begin{proof}
Let $D_{p,q}$ be the set of nw-diagrams of $T'_{p,q}$. 
We shall proceed by induction on $q$. If $q=1$, then $T'_{p,1}$ is
$$
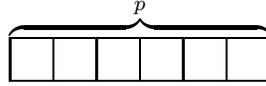
\overbrace{
\vbox{\offinterlineskip
\halign{#&#&#&#&#&#&#\cr
\htrait&\htrait&\htrait&\htrait&
\htrait&\htrait\cr
\sboite{}&\sboite{}&\sboite{}&\sboite{}&
\sboite{}&\sboite{}&\vtrait\cr
\htrait&\htrait&\htrait&\htrait&\htrait&\htrait\cr
}
}
}^{p}
$$
so we have 
$$
\sharp D_{p,q}=\mathcal T'_{p,q}=p+1=\frac{(p+q+1)!(p-q+2)}{q!(p+2)!}.
$$ 
Assume that $q>1$ and the claim is true for 
$q-1$. For $1\leqslant k\leqslant p-q+1$, let 
$$
S_{k}=\{ S\in D_{p,q}; (q,k)\in S \mbox{ and } (q,k+1)\not\in S\} 
$$
Then, $\sharp S_{k}=\mathcal T'_{p-k,q-1}$ and 
$L=\displaystyle\bigcup_{k=1}^{p-q+1}S_{k}$ is the set of
nw-diagrams containing at least a box in the last row of $T'_{p,q}$. Since 
$D_{p,q}$ is the disjoint union of $D_{p,q-1}$ and $L$, we obtain
that :
$$
\begin{array}{rcl}
\mathcal T'_{p,q}&=&\mathcal T'_{p,q-1}+ \sharp L = 
\displaystyle\sum_{i=0}^{p-q+1} \mathcal T'_{p-i,q-1} \\
&=& \displaystyle\sum_{k=q-1}^{p} \mathcal T'_{k,q-1} 
= \displaystyle\sum_{k=q-1}^{p} \frac{(k+q)!(k-q+3)}{(q-1)!(k+2)!}\\
&=& \displaystyle\frac{(p+q+1)!(p-q+2)}{q!(p+2)!} 
\end{array}
$$
where the last equality is a simple induction on $p\geqslant q$.
\end{proof}

\begin{definition}\label{defT_pq}
Let $p\geqslant q$ be two positive integers and $1\leqslant
l_1< l_2< \cdots< l_s \leqslant q+1$ be some other integers. 
Denote by $T_{p,q}(l_1,l_2,\dots,l_s)$ 
the new diagram obtained by adding
to $T_{p,q}$ the  boxes $(l_i,l_i-1)$, for $1\leqslant i \leqslant s$. 
For example, $T_{5,4}(2,4)$ is:
$$
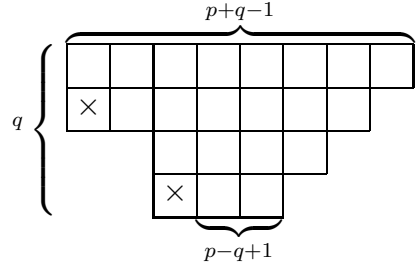
\vbox{\offinterlineskip
\halign{#\cr
${}^{q}\left\{ \vrule height36pt width0pt\right.$\cr
\vrule height1pt width0pt\cr
}
}
\overbrace{
\vbox{\offinterlineskip
\halign{#&#&#&#&#&#&#&#&#\cr
\htrait&\htrait&\htrait&\htrait&%
\htrait&\htrait&\htrait&\htrait&\cr
\sboite{}&\sboite{}&\sboite{}&\sboite{}&%
\sboite{}&\sboite{}&\sboite{}&\sboite{}&\vtrait\cr
\htrait&\htrait&\htrait&\htrait&%
\htrait&\htrait&\htrait&\htrait&\cr
\sboite{$\times$}&\sboite{}&\sboite{}&\sboite{}&\sboite{}&
\sboite{}&\sboite{}&\sboite{}&\cr
\htrait&\htrait&\htrait&\htrait&\htrait&\htrait&\htrait\cr
&&\sboite{}&\sboite{}&\sboite{}&\sboite{}&\vtrait\cr
&&\htrait&\htrait&\htrait&\htrait\cr
&&\sboite{$\times$}&\sboite{}&\sboite{}&\vtrait\cr
&&\multispan3\hrulefill&\cr
&&&\multispan2$\underbrace{\hboite{}\hboite{}}_{p-q+1}$&\cr
}
}
}^{p+q-1}
$$
where the added boxes are marked with a $\times$.
\end{definition}

\begin{prop}\label{compteB}
Let $p\geqslant q$ be two positive integers and $1\leqslant 
l_1< l_2< \cdots <l_s \leqslant q+1$ 
be some other integers, then the number of nw-diagrams of 
$T_{p,q}(l_1,l_2,\dots,l_s)$ is
$$
\left(\begin{array}{c}p+q \\ p\end{array}\right)
+\sum_{j=1}^s \mathcal T'_{p+q-l_j,l_j-1}.
$$
\end{prop}

\begin{proof}
Let $D_{p,q}(l_1,\dots,l_s)$ be the set of nw-diagrams of $T_{p,q}(l_1,\dots,l_s)$
and $\mathcal D_{p,q}(l_1,\dots,l_s)$ be its cardinality. 
Let $b_s=(l_s,l_s-1)$. Set 
$$
\begin{array}{l}
E=\{ S\in D_{p,q}(l_1,\dots,l_s); b_s\not\in S\}, \\
F=\{ S\in D_{p,q}(l_1,\dots,l_s); b_s\in S \mbox{ and } S\setminus \{b_s\}\in E\}, \\
G=\{ S\in D_{p,q}(l_1,\dots,l_s); b_s\in S \mbox{ and } S\setminus
\{b_s\} \not\in E\}. \\
\end{array}
$$
Then, we have clearly $\mathcal D_{p,q}(l_1,\dots,l_s)=\sharp E+\sharp F +\sharp G$.

If $S\in F$, then $S$ contains all the boxes north-west of $b_s$ and
the other boxes of $S$ are strictly north-east of $b_s$, so there
exists a bijection between $F$ and the set of nw-diagrams of
$T'_{p+q-l_s,l_s-1}$. For the example in definition \ref{defT_pq}, if
$S \in F$, $S$ is a nw-diagram of :
$$
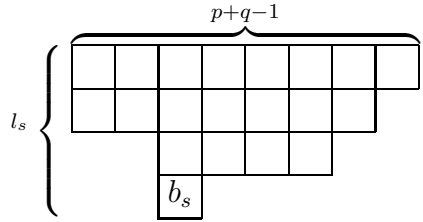
\vbox{\offinterlineskip
\halign{#\cr
${}^{l_s}\left\{ \vrule height36pt width0pt\right.$\cr
\vrule height1pt width0pt\cr
}
}
\overbrace{
\vbox{\offinterlineskip
\halign{#&#&#&#&#&#&#&#&#\cr
\htrait&\htrait&\htrait&\htrait&
\htrait&\htrait&\htrait&\htrait&\cr
\sboite{}&\sboite{}&\sboite{}&\sboite{}&
\sboite{}&\sboite{}&\sboite{}&\sboite{}&\vtrait\cr
\htrait&\htrait&\htrait&\htrait&
\htrait&\htrait&\htrait&\htrait&\cr
\sboite{}&\sboite{}&\sboite{}&\sboite{}&\sboite{}&
\sboite{}&\sboite{}&\sboite{}&\cr
\htrait&\htrait&\htrait&\htrait&\htrait&\htrait&\htrait\cr
&&\sboite{}&\sboite{}&\sboite{}&\sboite{}&\vtrait\cr
&&\htrait&\htrait&\htrait&\htrait\cr
&&\sboite{$b_s$}&\vtrait\cr
&&\htrait\cr
}
}
}^{p+q-1}
$$
containing $b_s$. Hence it suffices to count the nw-diagrams of the subdiagram
strictly north-east of $b_s$ :
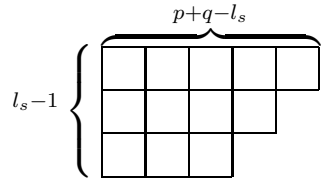
$$
\vbox{\offinterlineskip
\halign{#\cr
${}^{l_s-1}\left\{ \vrule height27pt width0pt\right.$\cr
\vrule height0pt width0pt\cr
}
}
\overbrace{
\vbox{\offinterlineskip
\halign{#&#&#&#&#&#\cr
\htrait&
\htrait&\htrait&\htrait&\htrait&\cr
\sboite{}&
\sboite{}&\sboite{}&\sboite{}&\sboite{}&\vtrait\cr
\htrait&
\htrait&\htrait&\htrait&\htrait&\cr
\sboite{}&\sboite{}&
\sboite{}&\sboite{}&\sboite{}&\cr
\htrait&\htrait&\htrait&\htrait\cr
\sboite{}&\sboite{}&\sboite{}&\vtrait\cr
\htrait&\htrait&\htrait\cr
}
}
}^{p+q-l_s}
$$
So by proposition \ref{compteA}, the cardinality of $F$ is 
$\mathcal T'_{p+q-l_s,l_s-1}$. 

If $S\in G$, then $S\setminus\{b_s\}$ is a nw-diagram of $T$ where $T=T_{p,q}$ if
$s=1$ and  $T=T_{p,q}(l_1,\dots,l_{s-1})$ if $s>1$. So the cardinality
of $G$ is the cardinality of the set of nw-diagrams in $T$ minus the cardinality of
the set $H$ of nw-diagrams having at most $l_s-1$ rows. Observe that
the elements of $H$ correspond to those of $E$. Hence, (by \cite{Pr}) 
$$
\sharp G=\left\{
\begin{array}{ll}
\mathcal D_{p,q}(l_1,\dots,l_{s-1}) -\sharp E &\mbox{if }
s>1, \\[2pt]
\displaystyle \left(
\begin{array}{c}p+q\\p\end{array}\right) -\sharp E &\mbox{if }
s=1.
\end{array}
\right.
$$
The result now follows easily by induction on $s$.
\end{proof}

\begin{notations}\label{Inotations}
Fix $I\subset \Pi$. Let
$I_1,\dots,I_s$ be the connected components of $I$ of 
cardinality $r_1, \dots,r_s$ respectively. For each connected component 
$I_j$, set $m_j=\min\{i;\alpha_i\in I_j\}$. Without loss of
generality, we shall assume that $m_1<m_2<\dots
<m_s$.
\end{notations}

If $\glie$ is of type $B_l$, then $T_{B_l}$ is $T_{l,l}$ and the boxes 
$(i,j)$ of $T_{B_l}$ are filled by 
the positive roots $t_{i,j}$, where
$$
t_{i,j}=\left\{
\begin{array}{ll}
\alpha_i+\cdots +2(\alpha_{j+1}+\cdots+\alpha_l), &1\leqslant j\leqslant l-1, \\
\alpha_i+\cdots+\alpha_{2l-j},&l\leqslant j\leqslant 2l-1.
\end{array}
\right.
$$
As before, for $I\subset \Pi$, we delete the boxes containing
elements of $\Delta_I$. For $j=1,\dots,s$, set 
\begin{equation}\label{l_jB}
l_j=m_j-\sum_{k=1}^{j-1} r_k.
\end{equation}
Regroup the equivalent classes of $\sim_I$ proceeding simple 
root by simple root : for each 
$\alpha_i\in I\setminus\{\alpha_l\}$, we first regroup rows $i$ and
$i+1$ and if $i\not=1$, we regroup column $2l-i$ 
and column $2l-i+1$, then the columns $i-1$ and $i$. 
If $\alpha_l\in I$, we also regroup the columns $l-1$, $l$ and $l+1$. We obtain that 
$T_{B_l}^I$ is a  diagram 
of shape $T_{l-\sharp I,l-\sharp I}(l_1,\dots,l_n)$, 
where the $l_i$ are defined as above and, 
$n=s-1$ if $\alpha_l\in I$ and $n=s$ if $\alpha_l\not\in I$. For
example, for $B_5$ and $I=\{\alpha_2,\alpha_3,\alpha_5\}$, we have :
$$
\begin{array}{ccc}
\vbox{\offinterlineskip
\halign{#&#&#&#&#&#&#&#&#&#&#&#&#&#&#&#&#&#&#\cr
\htrait&\htrait&\htrait&\htrait&\htrait&
\htrait&\htrait&\htrait&\htrait\cr
\sboite{$t_{1,1}$}
&\hboite{$t_{1,2}$}
&\hboite{$t_{1,3}$}
&\sboite{$t_{1,4}$}&
\hboite{$t_{1,5}$}
&\hboite{$t_{1,6}$}
&\sboite{$t_{1,7}$}
&\hboite{$t_{1,8}$}
&\hboite{$t_{1,9}$}&\vtrait\cr
\htrait&\htrait&\htrait&
\htrait&\htrait&\htrait&\htrait&\htrait&\htrait\cr
&\sboite{$t_{2,2}$}
&\hboite{$t_{2,3}$}
&\sboite{$t_{2,4}$}
&\hboite{$t_{2,5}$}
&\hboite{$t_{2,6}$}
&\vtrait\cr
&\htrait\cr
&&\sboite{$t_{3,3}$}&\sboite{$t_{3,4}$}&\hboite{$t_{3,5}$}&\hboite{$t_{3,6}$}&\vtrait\cr
&&\htrait\cr
&&&\sboite{$t_{4,4}$}&\hboite{$t_{4,5}$}&\hboite{$t_{4,6}$}&\vtrait\cr
&&&\htrait&\htrait&\htrait\cr
}
}
&
\vbox{\offinterlineskip
\halign{#\cr
$\longleftrightarrow$\cr
\vrule height30pt width0pt\cr
}
}
&\vbox{\offinterlineskip
\halign{#&#&#&#\cr
\htrait&\htrait&\htrait\cr
\sboite{}&\sboite{}&\sboite{}&\vtrait\cr
\htrait&\htrait&\htrait\cr
\sboite{$b_1$}&\sboite{}&\vtrait\cr
\htrait&\htrait\cr
\vrule height16pt  width0pt\cr
}
}
\end{array}
$$

It follows from proposition \ref{compteB} that :

\begin{prop}\label{enumerationB}
Let $I\subset\Pi$ be of cardinality $r$. 
Let $\mathcal S_{B_l}^I$ be the set of all nw-diagrams of
$T_{B_l}^I$. Then, the cardinality of $\mathcal S_{B_l}^I$ is 
$$
(l-r+1)\mathcal C_{l-r}+
\sum_{j=1}^{n} \mathcal T'_{2(l-r)-l_j,l_j-1}
$$
where $n=s-1$ if $\alpha_l\in I$, and $n=s$ otherwise.
\end{prop}

If $\glie$ is of type $D_l$, then $T_{D_l}$ is $T_{l,l-1}$, and the boxes 
$(i,j)$ of $T_{D_l}$ are filled by 
the positive roots $t_{i,j}$, where
$$
t_{i,j}=\left\{
\begin{array}{ll}
\alpha_i+\cdots +2(\alpha_{j+1}+\cdots+\alpha_{l-2})
+\alpha_{l-1}+\alpha_l, &1\leqslant j\leqslant l-2, \\
\alpha_i+\cdots+\alpha_{l-2}+\alpha_l, &j=l-1, \\
\alpha_i+\cdots+\alpha_{2l-j}, &l\leqslant j\leqslant 2l-1.
\end{array}
\right.
$$
For $I\subset \Pi$, we first delete the boxes containing elements of
$\Delta_I$. For $j=1,\dots,s$, set 
\begin{equation}\label{l_jD}
l_j=\left\{
\begin{array}{ll}
\displaystyle m_j-\sum_{k=1}^{j-1} r_k &\mbox{if }
j\not=s \mbox{ or }I_s\not=\{\alpha_{l}\},\\
\displaystyle m_j-\sum_{k=1}^{j-1} r_k-1 &\mbox{if } 
j=s \mbox{ and } I_s=\{\alpha_{l}\}.
\end{array}
\right.
\end{equation}
Regroup the equivalent classes of $\sim_I$ proceeding simple 
root by simple root : for each 
$\alpha_i\in I\setminus\{\alpha_{l-1},\alpha_l\}$, we first regroup
the rows $i$ and $i+1$ and if $i\not=1$, we regroup column $2l-i-1$ 
and column $2l-i$, and then the columns $i-1$ and $i$. 

If $\alpha_{l-1}\in I$, but $\alpha_l\not\in I$, then we
regroup the columns $l-2,l-1$ and columns $l,l+1$.
 
If $\alpha_{l}\in I$, but $\alpha_{l-1}\not\in I$, then we first reverse the
columns $l-1$ and $l$, and then we regroup the (new) columns $l-2,l-1$
and columns $l,l+1$. 

If $\{\alpha_{l-1},\alpha_l\}\subset I$, then we
regroup the four columns $l-2,l-1$, $l$ and $l+1$. 

We obtain that if $\{\alpha_{l-1},\alpha_l\}\not\subset I$, then 
$T_{D_l}^I$ is a diagram  of shape 
$T_{l-\sharp I,l-\sharp I-1}(l_1,\dots,l_s)$, where the $l_i$
are defined as above. 
If $\{\alpha_{l-1},\alpha_l\}\subset I$, then $T_{D_l}^I$ is a
diagram of shape $T_{l-\sharp I,l-\sharp I}(l_1,\dots,l_{s-1})$.

In the following examples, we denote by $i$ the simple root $\alpha_i$ and
by $i^2$ the element $2\alpha_i$. We consider, $X=D_5$ and $I$ is
respectively $\{\alpha_1,\alpha_2,\alpha_5\}$ and
$\{\alpha_2,\alpha_4,\alpha_5\}$ :

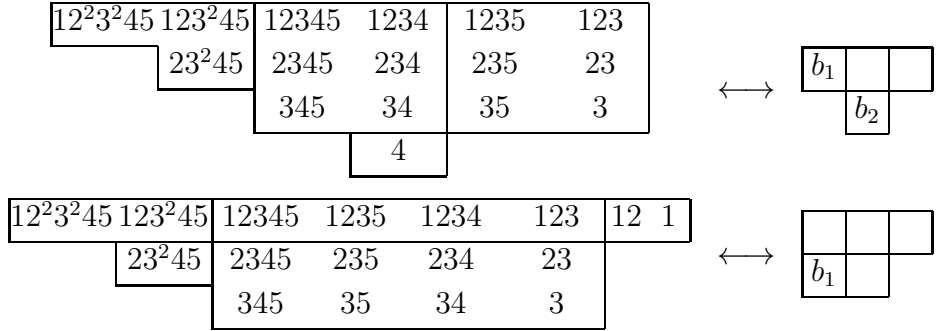
$$
\begin{array}{ccc}
\vbox{\offinterlineskip
\halign{#&#&#&#&#&#&#\cr
\htrait&\htrait&\htrait&\htrait&\htrait&\htrait\cr
\sb{$12^23^245$}
&\hb{$123^245$}&\sb{$12345$}&\hb{$1234$}&\sb{$1235$}
&\hb{$123$}&\vt\cr
\htrait&&&&&\cr
&\sb{$23^245$}&\sb{$2345$}&\hb{$234$}&\sb{$235$}
&\hb{$23$}&\vt\cr
&\htrait&&&&\cr
&&\sb{$345$}&\hb{$34$}&\sb{$35$}
&\hb{$3$}&\vt\cr
&&\htrait&\htrait&\htrait&\htrait\cr
&&&\sb{$4$}&\vt\cr
&&&\htrait\cr
}
} 
&
\vbox{\offinterlineskip
\halign{#\cr
$\longleftrightarrow$\cr
\vrule height30pt width0pt\cr
}
}
&\vbox{\offinterlineskip
\halign{#&#&#&#\cr
\htrait&\htrait&\htrait\cr
\sboite{$b_1$}&\sboite{}&\sboite{}&\vtrait\cr
\htrait&\htrait&\htrait\cr
&\sboite{$b_2$}&\vtrait\cr
&\htrait\cr
\vrule height16pt  width0pt\cr
}
} \\[4pt]
\vbox{\offinterlineskip
\halign{#&#&#&#&#&#&#&#&#\cr
\htrait&\htrait&\htrait&\htrait&\htrait&\htrait&\htrait&\htrait\cr
\sb{$12^23^245$}
&\hb{$123^245$}&\sb{$12345$}&\hb{$1235$}&\hb{$1234$}
&\hb{$123$}&\sboite{$12$}&\hboite{$1$}&\vt\cr
\htrait&\htrait&\htrait&\htrait&\htrait&\htrait&\htrait&\htrait\cr
&\sb{$23^245$}&\sb{$2345$}&\hb{$235$}&\hb{$234$}
&\hb{$23$}&\vt\cr
&\htrait\cr
&&\sb{$345$}&\hb{$35$}&\hb{$34$}
&\hb{$3$}&\vt\cr
&&\htrait&\htrait&\htrait&\htrait\cr
}
} 
&
\vbox{\offinterlineskip
\halign{#\cr
$\longleftrightarrow$\cr
\vrule height25pt width0pt\cr
}
}
&\vbox{\offinterlineskip
\halign{#&#&#&#\cr
\htrait&\htrait&\htrait\cr
\sboite{}&\sboite{}&\sboite{}&\vtrait\cr
\htrait&\htrait&\htrait\cr
\sboite{$b_1$}&\sboite{}&\vt\cr
\htrait&\htrait\cr
\vrule height12pt  width0pt\cr
}
}

\end{array}
$$

\begin{definition}
For a subdiagram $L$ of $ T_{D_l}^I$, we shall denote by 
$L^{\bullet}$ the set of boxes of $L$
obtained from $L$ by exchanging columns $l-r-1$ and $l-r$ (resp. $l-r$
and $l-r+1$) if $\alpha_1\not\in I$ (resp. if $\alpha_1\in I$). 

If $L^{\bullet}$ is a nw-diagram of $T_{D_l}^I$, then we say that $L$ is a
$\bullet$-nw-diagram of $T_{D_l}^I$.
\end{definition}

\begin{prop}\label{enumerationD}
Let $I\subset\Pi$ be of cardinality $r$. Let 
$\mathcal S_{D_l}^I$ be the set of nw-diagrams of
$T_{D_l}^I$ if $\{\alpha_{l-1},\alpha_l\}\cap I\not=\emptyset$, and be
the union of the set of nw-diagrams and the set of
$\bullet$-nw-diagrams of 
$T_{D_l}^I$ if 
$\{\alpha_{l-1},\alpha_l\}\cap I=\emptyset$. Then, the cardinality of 
$\mathcal S_{D_l}^I$ is 

\begin{enumerate}[{\rm(i)}]

\item $\displaystyle (3(l-r)-2)\mathcal C_{l-r-1}  +
\displaystyle\sum_{j=1}^{s} \mathcal T'_{2(l-r)-l_j-1,l_j-2}+
\mathcal T'_{2(l-r)-l_j-1,l_j-1}$, if $t=0$,

\item $\displaystyle \frac{l-r+1}{2} \mathcal C_{l-r} 
+
\displaystyle\sum_{j=1}^{s} \mathcal T'_{2(l-r)-l_j-1,l_j-1}$, if
$t=1$,

\item $(l-r+1)\mathcal C_{l-r}
+\displaystyle
\sum_{j=1}^{s-1} \mathcal T'_{2(l-r)-l_j,l_j-1}$, if
$t=2$,
\end{enumerate}
where $t=\sharp(\{\alpha_{l-1},\alpha_l\}\cap I)$.
\end{prop}

\begin{proof}
Assume first that $\{\alpha_{l-1},\alpha_l\}\cap I=\emptyset$. Note that,
the elements of $\calf_I$ are in bijection with the subdiagrams $S$ of 
$T_{D_l}^I=T_{l-r,l-r-1}(l_1,\dots,l_s)$ such that
either $S$ or $S^{\bullet}$ is a nw-diagram. Let 
$$
\begin{array}{l}
E_1=\mbox{ the set of nw-diagrams of } T_{D_l}^I, \\
E_2=(\mbox{the set of  $\bullet$-nw-diagrams of }T_{D_l}^I)\setminus E_1.
\end{array}
$$
So $S_{D_l}^I=E_1\cup E_2$ (disjoint union).
By proposition \ref{compteB}, we have :
$$
\sharp E_1=\displaystyle \frac{l-r+1}{2}\mathcal C_{l-r}+
\sum_{j=1}^{s}\mathcal T'_{2(l-r)-l_j-1,l_j-1}.
$$
On the other hand, the number of elements of $E_2$  is 
$\sharp E_1-\sharp F$, where $F$ is the set of elements of $E_1$ having columns 
$l-r-1$ and $l-r$ (resp. $l-r$ and $l-r+1$) of the same length if
$\alpha_1\not\in I$ (resp. if $\alpha_1\in I$). 

Clearly, the number of elements of $F$ 
is exactly the number of nw-diagrams of the diagram  obtained from
$T_{D_l}^I$ by removing the $(l-r)$-th (resp. $(l-r+1)$-th) column if 
$\alpha_1\not\in I$ (resp. if $\alpha_1\in I$). So, by proposition \ref{compteB},
$$
\sharp F=
\displaystyle (l-r)\mathcal C_{l-r-1}
+
\sum_{j=1}^{s} \mathcal T'_{2(l-r)-l_j-2,l_j-1}. 
$$
We obtain therefore the result since we have the equality : 
$$
\begin{array}{c}
\mathcal T'_{2(l-r)-l_j-1,l_j-1}-\mathcal T'_{2(l-r)-l_j-2,l_j-1}
=\mathcal T'_{2(l-r)-l_j-1,l_j-2}.
\end{array}
$$

If $\alpha_{l-1}$ or $\alpha_l\in I$, then there is no column
reversing. Then the result follows from proposition \ref{compteB}
according to the shape of $T_{D_l}^I$.
\end{proof}

As in \cite{CP}, we have clearly a bijection between $\calf_I$ and $\mathcal
S_X^I$. It follows from propositions \ref{enumerationA},
\ref{enumerationC}, \ref{enumerationB}, \ref{enumerationD}, that we have :

\begin{theoreme}\label{thenumeration}
Let $I \subset \Pi$ of cardinality $r$, $X$ be the type of $\glie$ and $s,l_j$ as defined
in \ref{Inotations}, \eqref{l_jB} and \eqref{l_jD}.

If $X=A_l$, then 
$$
\sharp \calf_I=
\mathcal C_{l-r+1}.
$$

If $X=B_l$, then 
$$
\sharp \calf_I=(l-r+1)\mathcal C_{l-r}+
\sum_{j=1}^{n} \mathcal T'_{2(l-r)-l_j,l_j-1},
$$
where $n=s-1$ if $ \alpha_l\in I$, and $n=s$ otherwise.

If $X=C_l$, then 
$$
\sharp \calf_I=
\left\{
\begin{array}{ll}
(l-r +1)\mathcal C_{l-r} &\mbox{if }\  \alpha_l\not\in
  I, \\[2pt]
\displaystyle \frac{l-r +2}{2}\mathcal C_{l-r+1} 
&\mbox{if }\  \alpha_l\in I.
\end{array}
\right.
$$

If $X=D_l$, then 
$$
\sharp \calf_I=
\left\{
\begin{array}{ll}
\displaystyle \frac{l-r+1}{2} \mathcal C_{l-r} 
+
\displaystyle\sum_{j=1}^{s} \mathcal T'_{2(l-r)-l_j-1,l_j-1}
&
\mbox{if }  \sharp\{\alpha_{l-1},\alpha_l\}\cap I=1\\
 (l-r+1)\mathcal C_{l-r}
+\displaystyle
\sum_{j=1}^{s-1} \mathcal T'_{2(l-r)-l_j,l_j-1}
&
\mbox{if } \{\alpha_{l-1},\alpha_l\}\subset I,\\
\displaystyle (3(l-r)-2)\mathcal C_{l-r-1} 
+&\\
\displaystyle\sum_{j=1}^{s} \mathcal T'_{2(l-r)-l_j-1,l_j-2}+
\mathcal T'_{2(l-r)-l_j-1,l_j-1}
&\mbox{otherwise.}
\end{array}
\right.
$$
\end{theoreme}

%--------------------------------------------------------------------------------------------------------------------------
\subsection{Abelian ideals}

We have already determined in theorem \ref{theoremeab} the number of 
abelian ideals for type $A$ and $C$. We shall now enumerate the
abelian ideals of $\plie_I$ using diagrams when $\glie$ is of type
$B$ or $D$. Observe that a similar argument could be used to enumerate
abelian ideals in type $A$ and $C$.

\begin{definition}
Let $p$ be a positive integer and $R_{p}$ be the diagram of shape
$[p,p-1,\dots,1]$ arranged in
the  following way :
$$
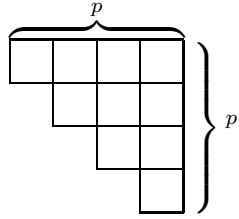
\overbrace{
\vbox{\offinterlineskip
\halign{#&#&#&#&#\cr
\htrait&\htrait&\htrait&\htrait\cr
\sboite{}&\sboite{}&\sboite{}&\sboite{}&\vtrait\cr
\htrait&\htrait&\htrait&\htrait\cr
&\sboite{}&\sboite{}&\sboite{}&\vtrait\cr
&\htrait&\htrait&\htrait\cr
&&\sboite{}&\sboite{}&\vtrait\cr
&&\htrait&\htrait\cr
&&&\sboite{}&\vtrait\cr
&&&\htrait\cr
}
}
}^{p}
\vbox{\offinterlineskip
\halign{#\cr
$\left. \vrule height36pt width0pt\right\}{}^{p}$\cr
\vrule height0pt width0pt\cr
}
}
$$
\end{definition}

\begin{prop}\label{R_p}
The number of nw-diagrams of $R_{p}$ is $2^p$.
\end{prop}

\begin{proof}
We shall procced by induction on $p$. If $p=1$, the result is
clear. Assume that $p>1$ and the claim is true for $p-1$. Let $b$ be
the box $(1,p)$ and 
\begin{eqnarray*}
E&=&\mbox{the set of nw-diagrams of $R_p$ which do not contain $b$, } \\
F&=&\mbox{the set of nw-diagrams of $R_p$ which contain $b$. }
\end{eqnarray*}
Then, the number of nw-diagrams of $R_p$ is $\sharp E+\sharp
F$. Furthermore, by the induction hypothesis, we have $\sharp
E=2^{p-1}=\sharp F$, and we obtain the result.
\end{proof}

\begin{definition}\label{defR}
Let $p$ be a positive integer and $1\leqslant
l_1< l_2< \cdots< l_s \leqslant p+1$ 
be some other integers. Denote by $R_{p}(l_1,l_2,\dots,l_s)$ 
the new diagram obtained by adding
to $R_{p}$ the  boxes $(l_i,l_i-1)$, for $1\leqslant i \leqslant s$. 
For example, $R_{4}(3)$ :
$$
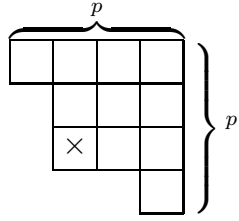
\overbrace{
\vbox{\offinterlineskip
\halign{#&#&#&#&#\cr
\htrait&\htrait&\htrait&\htrait\cr
\sboite{}&\sboite{}&\sboite{}&\sboite{}&\vtrait\cr
\htrait&\htrait&\htrait&\htrait\cr
&\sboite{}&\sboite{}&\sboite{}&\vtrait\cr
&\htrait&\htrait&\htrait\cr
&\sboite{$\times$}&\sboite{}&\sboite{}&\vtrait\cr
&\htrait&\htrait&\htrait\cr
&&&\sboite{}&\vtrait\cr
&&&\htrait\cr
}
}
}^{p}
\vbox{\offinterlineskip
\halign{#\cr
$\left. \vrule height36pt width0pt\right\}{}^{p}$\cr
\vrule height0pt width0pt\cr
}
}
$$

\end{definition}

\begin{prop}\label{compteab}
Let $p$ be a positive integer and $1\leqslant 
l_1< l_2< \cdots <l_s \leqslant p+1$ be
some other integers, then the number of nw-diagrams of 
$R_{p}(l_1,l_2,\dots,l_s)$ is
$$
2^{p}+\sum_{j=1}^s \left(\begin{array}{c}p \\ l_j-1\end{array}\right).
$$
\end{prop}

\begin{proof}
Let $D_{p}(l_1,\dots,l_s)$ be the set of nw-diagrams of
$R_{p}(l_1,\dots,l_s)$ and $\mathcal D_{p}(l_1,\dots,l_s)$ be its
cardinality. 
Let $b_s=(l_s,l_s-1)$. Set 
$$
\begin{array}{l}
E=\{ S\in D_{p}(l_1,\dots,l_s); b_s\not\in S\} \\
F=\{ S\in D_{p}(l_1,\dots,l_s); b_s\in S \mbox{ and } S\setminus \{b_s\}\in E\} \\
G=\{ S\in D_{p}(l_1,\dots,l_s); b_s\in S \mbox{ and } S\setminus
\{b_s\} \not\in E\} \\
\end{array}
$$
Then we have clearly $\mathcal D_{p}(l_1,\dots,l_s)=\sharp E+\sharp
F+\sharp G$.

If $S\in F$, then $S$ contains all the boxes north-west of $b_s$ and
the other boxes of $S$ are strictly north-east of $b_s$, so there
exists a bijection between $F$ and the set of nw-diagrams of $T$ where
$T$ is a diagram whose  shape is a rectangle containing 
$p-l_j+1$ columns and $l_j-1$ rows. For the example in definition
\ref{defR}, if $S \in F$, $S$ is a nw-diagram of :

$$
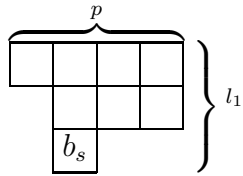
\overbrace{
\vbox{\offinterlineskip
\halign{#&#&#&#&#\cr
\htrait&\htrait&\htrait&\htrait\cr
\sboite{}&\sboite{}&\sboite{}&\sboite{}&\vtrait\cr
\htrait&\htrait&\htrait&\htrait\cr
&\sboite{}&\sboite{}&\sboite{}&\vtrait\cr
&\htrait&\htrait&\htrait\cr
&\sboite{$b_s$}&\vtrait&&\cr
&\htrait\cr
}
}
}^{p}
\vbox{\offinterlineskip
\halign{#\cr
$\left. \vrule height27pt width0pt\right\}{}^{l_1}$\cr
\vrule height0pt width0pt\cr
}
}
$$
containing $b_s$. Hence it suffices to count the nw-diagrams of the
rectangular subdiagram strictly north-east of $b_s$ :
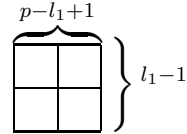
$$
\overbrace{
\vbox{\offinterlineskip
\halign{#&#&#\cr
\htrait&\htrait\cr
\sboite{}&\sboite{}&\vtrait\cr
\htrait&\htrait\cr
\sboite{}&\sboite{}&\vtrait\cr
\htrait&\htrait\cr
}
}
}^{p-l_1+1}
\vbox{\offinterlineskip
\halign{#\cr
$\left. \vrule height19pt width0pt\right\}{}^{l_1-1}$\cr
\vrule height0pt width0pt\cr
}
}
$$
So by \cite{Pr} the cardinality of $F$ is 
$\left(\begin{array}{c}p \\ l_j-1\end{array}\right)$. 

If $S\in G$, then $S\setminus\{b_s\}$ is a nw-diagram of $L$ where $L=R_{p}$ if
$s=1$ and  $L=R_{p}(l_1,\dots,l_{s-1})$ if $s>1$. So the cardinality
of $G$ is the cardinality of the set of nw-diagrams in $L$ minus the cardinality of
the set $H$ of nw-diagrams having at most $l_s-1$ rows. Observe that
the elements of $H$ correspond to those of $E$. Hence, by proposition \ref{R_p} 
$$
\sharp G=\left\{
\begin{array}{ll}
\mathcal D_{p,q}(l_1,\dots,l_{s-1}) -\sharp E &\mbox{if }
s>1, \\[2pt]
\displaystyle 2^p-\sharp E &\mbox{if }
s=1.
\end{array}
\right.
$$
The result now follows easily by induction on $s$.
\end{proof}

Let $\calf_I^{ab}=\{\Phi\in \calf_I; \ilie_{\Phi} \mbox{ is abelian
}\}.$ If $S$ is a subdiagram of a diagram, let 
$$
\tau_h^{S}=\max\{k;(h,k)\in S\},
$$
so $(h,\tau_h^{S})$ is the right most box in the $h$-th row of $S$.

\begin{prop}\label{condabelienB}
Assume that $\glie$ is of type $B_l$. Let $I \subset \Pi$ be of cardinality
$r$. Consider 
$\Phi \in \calf_I$ and $S$ its corresponding nw-diagram in $T_{B_l}^I$.
Then $\Phi\in \calf_I^{ab}$ if and only if 

{\rm (a)} $\tau_1^{S} \leq l-r$ 
if $\alpha_1\in I$,

{\rm (b)} $\tau_1^{S} + \tau_2^{S} \leq 2(l-r)-1$ 
if $\alpha_1\not\in I$.
\end{prop}

\begin{proof}
Let $S_0$ be the corresponding nw-diagram of $\Phi$ in $T_{B_l}^{\emptyset}$, then
by \cite{CP}, we have $\Phi\in \calf_{\emptyset}^{ab}$ if and only if 
$\tau_1^{S_0} +\tau_2^{S_0}\leq 2l-1$.

If $\alpha_1\in I$, then $\tau_1^{S_0}=\tau_2^{S_0}$. The regrouping
process reduces the number of columns on the left of column $l$ of
$T_{B_l}^{\emptyset}$ by one for each simple root in
$I\setminus\{\alpha_1\}$. It follows that $\Phi\in \calf_I^{ab}$ if
and only if $\tau_1^{S} \leq l-r$.

The argument is similar for the case $\alpha_1\not\in I$.
\end{proof}

\begin{prop}\label{propenumerationB}
Assume that $\glie$ is of type
$B_l$. Let $I\subset \Pi$ be of cardinality $r$, and $l_1,\dots,l_s$ be as
defined in \eqref{l_jB}. Then we have :
$$
\sharp \calf_I^{ab}=
\left\{
\begin{array}{ll}
\displaystyle 2^{l-r}+\sum_{j=1}^n 2
\left(\begin{array}{c} l-r-1 \\l_j-1\end{array}\right)
&\mbox{if } \alpha_1\not\in I, \\
\displaystyle 2^{l-r-1}+\sum_{j=1}^n 
\left(\begin{array}{c} l-r-1 \\l_j-1\end{array}\right)
&\mbox{if } \alpha_1\in I, 
\end{array}
\right.
$$
where $n=s$ if $\alpha_l\not\in I$ and $n=s-1$ if $\alpha_l\in I$.
\end{prop}

\begin{proof} 
Recall that $T_{B_l}^I$ is of shape $T_{l-r,l-r}(l_1,\dots,l_n)$,
where $n=s$ if $\alpha_l\not\in I$ and $n=s-1$ if $\alpha_l\in I$.

Let $\Phi\in \calf_I$ and $S$ be the nw-diagram of
$T_{B_l}^{I}$ corresponding to $\Phi$.

Assume that $\alpha_1\in I$, then $l_1=1$. By proposition \ref{condabelienB}, $S$ is
in the left hand half of $T_{B_l}^{I}$, so it is a nw-diagram of
$R_{l-r-1}(1,\dots,l_n)$. We then obtain the result by proposition \ref{compteab}.

Assume that $\alpha_1\not\in I$. Let $E$ be the set of nw-diagrams of
$T_{B_l}^{I}$ associated to elements of $\calf_I^{ab}$. Set 
$$
\begin{array}{l}
P=\{S\in E; \tau_1^{S}\leqslant l-r-1\}, \\
Q=\{S\in E; \tau_1^{S}> l-r-1\}.
\end{array}
$$
Then, we have $\sharp E=\sharp P+\sharp Q$.

If $S\in P$, then $S$ is included in the left hand half of
$T_{B_l}^{I}$, so 
$$
\sharp P=2^{l-r-1}+\sum_{j=1}^n 
\left(\begin{array}{c} l-r-1 \\l_j-1\end{array}\right)
$$ 
by proposition \ref{compteab}.

For $i=l-r,\dots,2(l-r)-1$, let $Q_i=\{S\in Q; \tau_1^S=i\}$ and $P_i=\{S\in P;
\tau_1^S=2(l-r)-1-i\}$. We then have : 
$$
\displaystyle
Q=\bigcup_{i=l-r}^{2(l-r)-1} Q_i \ \mbox{ and }\ 
P=\bigcup_{i=l-r}^{2(l-r)-1} P_i.
$$
For $i=l-r,\dots,2(l-r)-1$, we have an obvious bijection between $P_i$ and
$Q_i$ given by the adding or removing of boxes $(1,2(l-r)-i),
\dots,(1,i)$. Therefore $\sharp P=\sharp Q$ and the result follows.
\end{proof}

\begin{prop}\label{condabelienD}
Assume that $\glie$ is of type $D_l$. Let $I \subset \Pi$ be of cardinality
$r$. Consider $\Phi \subset \calf_I$ and $S_{\Phi}$ its corresponding subdiagram
in $T_{D_l}^I$. Set $S=S_{\Phi}$ if $S_{\Phi}$ is a nw-diagram and 
$S=S_{\Phi}^{\bullet}$ if $S_{\Phi}^{\bullet}$ is a nw-diagram. 
Then $\Phi\in \calf_I^{ab}$ if and only if 

{\rm (a)} $\tau_1^{S} \leq l-r$ if $\alpha_1\in I$,

{\rm (b)} $\tau_1^{S} + \tau_2^{S} \leq 2(l-r)-2$ 
if $\alpha_1\not\in I$. 
\end{prop}

\begin{proof}
If $I=\emptyset$, set $S=S_0$, then by \cite{CP}, we have 
$\Phi\in \calf_{\emptyset}^{ab}$ if and only if 
$\tau_1^{S_0} +\tau_2^{S_0}\leq 2l-2$.

Assume that $\alpha_1\in I$, then $\tau_1^{S_0}=\tau_2^{S_0}$. The regrouping
process reduces the number of columns of the left of column $l$ of
$T_{D_l}^{\emptyset}$ by one for each simple root in
$I\setminus\{\alpha_1\}$. It follows that $\Phi\in \calf_I^{ab}$ if
and only if $\tau_1^{S} \leq l-r$.

The argument is similar for the case $\alpha_1\not\in I$.
\end{proof}

\begin{prop}\label{propenumerationD}
Assume that $\glie$ is of type $D_l$. Let $I \subset \Pi$ be of cardinality
$r$ and $l_1,\dots,l_s$ be as defined in \eqref{l_jD}. 
Set $t=\sharp (\{\alpha_{l-1},\alpha_l\}\cap I)$. If $\alpha_1 \in
I$, then the cardinality of $\calf_I^{ab}$ is :

{\rm(i)}  $\displaystyle
2^{l-r}-2^{l-r-2}+\sum_{j=1}^{s}\left[2\left(\begin{array}{c} l-r-1
    \\l_j-1\end{array}\right)-
\left(\begin{array}{c} l-r-2 \\l_j-1\end{array}\right)\right]$, 
if $t=0$, 

{\rm(ii)}  $\displaystyle
2^{l-r-1}+\sum_{j=1}^{s} \left(\begin{array}{c} l-r-1
  \\l_j-1\end{array}\right)$, if $t=1$,

{\rm(iii)}  $\displaystyle
2^{l-r-1}+\sum_{j=1}^{s-1} \left(\begin{array}{c} l-r-1
  \\l_j-1\end{array}\right)$, if $t=2$.

\noindent
If $\alpha_1 \not\in I$, then the cardinality of $\calf_I^{ab}$ is :

{\rm(iv)}  $\displaystyle
2^{l-r}+\sum_{j=1}^{s}2\left(\begin{array}{c} l-r-1
    \\l_j-1\end{array}\right)$, if $t=0$, 

{\rm(v)}  $\displaystyle 2^{l-r-1}+2^{l-r-2}+\sum_{j=1}^{s} 
\left(\begin{array}{c} l-r-1 \\l_j-1\end{array}\right)
+ \sum_{j=1}^{s-1} 
\left(\begin{array}{c} l-r-2 \\l_j-1\end{array}\right)$, if 
$t=1$, 

{\rm(vi)}  $\displaystyle 2^{l-r}+2\sum_{j=1}^{s-1} \left(\begin{array}{c} l-r-1
  \\l_j-1\end{array}\right) $, if $t=2$.

\end{prop}

\begin{proof}
We proceed as in the case of type $B_l$ but here, we need to take
into account column reversing. 

Recall that if $t=0$ or $1$, $T_{D_l}^I$ is
of shape $T_{l-r,l-r-1}(l_1,\dots,l_s)$ and if $t=2$, $T_{D_l}^I$ is
of shape $T_{l-r,l-r}(l_1,\dots,l_{s-1})$.

Let $\mathcal S^{ab}_I$ be the set of subdiagrams of $T_{D_l}^I$
corresponding to  elements of $\calf_I^{ab}$. The shape of elements
of $\mathcal S^{ab}_I$ is
conditioned by proposition \ref{condabelienD}. Let 
$$
\begin{array}{l}
E_1=\mbox{the set of nw-diagrams in }\mathcal S^{ab}_I, \\
E_2=(\mbox{the set of $\bullet$-nw-diagrams in }\mathcal S^{ab}_I)\setminus E_1. \
\end{array}
$$
Consider $\Phi \in \calf_I^{ab}$ and $S$ its corresponding subdiagram
in $\mathcal S^{ab}_I$.

First assume that $\alpha_1\in I$, then $l_1=1$. If $S\in E_1$, by
proposition \ref{condabelienD}, $S$ is
in the left hand half of $T_{D_l}^{I}$, so it is a nw-diagram of
$R_{l-r-1}(1,\dots,l_n)$, where $n=s$ if $t=0,1$ and $n=s-1$ if $t=2$.
Hence, by proposition \ref{compteab}, we have 
$$
\sharp E_1=2^{l-r-1}+\sum_{j=1}^n \left(\begin{array}{c} l-r-1\\
  l_j-1\end{array}\right).
$$

If $t\not=0$, there is no column reversing, so $E_2=\emptyset$. If
$t=0$, the number of elements of $E_2$  is 
$\sharp E_1-\sharp (F\cap E_1)$, where $F$ 
is the set of nw-diagrams of $T_{D_l}^I$ having columns 
$l-r$ and $l-r+1$ of the same length. 

Clearly, the number of elements of $F$ 
is exactly the number of nw-diagrams of the diagram  obtained from
$T_{D_l}^I$ by removing the $(l-r+1)$-th column. So, by proposition
\ref{condabelienD}, the set of elements which are in $F\cap E_1$ is in
bijection with the set of nw-diagrams of $R_{l-r-2}(1,\dots,l_s).$ 
So by proposition
\ref{compteab}, we obtain :
$$
\sharp F=
2^{l-r-2}+\sum_{j=1}^{s}\left(\begin{array}{c}l-r-2\\
  l_j-1\end{array}\right).  
$$
We obtain therefore the result.

Now assume that $\alpha_1\not\in I$. Set
$$
\begin{array}{l}
P=\{S\in E_1; \tau_1^{S}\leqslant l-r-1\}, \\
\widetilde P=\{S\in E_1; \tau_1^{S}\leqslant l-r-2\}, \\
Q=\{S\in E_1; \tau_1^{S}> l-r-1\}.
\end{array}
$$
Then, we have $\sharp E_1=\sharp P+\sharp Q$.

First assume that $t=0$ or $1$. If $S\in P$, then $S$ is 
included in the left hand half of
$T_{D_l}^{I}$, so by proposition \ref{compteab}, we have 
$$\sharp P=2^{l-r-1}+\sum_{j=1}^s 
\left(\begin{array}{c} l-r-1 \\l_j-1\end{array}\right)
$$ 

For $i=l-r,\dots,2(l-r)-2$, let $Q_i=\{S\in Q; \tau_1^S=i\}$ and
$\widetilde P_i=\{S\in P;
\tau_1^S=2(l-r)-2-i\}$. We then have : 
$$
Q=\bigcup_{i=l-r}^{2(l-r)-2} Q_i \ \mbox{ and } \ 
\widetilde P=\bigcup_{i=l-r}^{2(l-r)-2} \widetilde P_i.
$$
For $i=l-r,\dots,2(l-r)-2$, we have an obvious bijection between
$\widetilde P_i$ and $Q_i$ given by the adding or removing of boxes $(1,2(l-r)-i),
\dots,(1,i)$. Therefore $\sharp \widetilde P=\sharp Q$ and by
proposition \ref{compteab}, we have 
$$
\sharp \widetilde P =2^{l-r-2}+\sum_{j=1}^s 
\left(\begin{array}{c} l-r-2 \\l_j-1\end{array}\right).
$$
If $t=1$, there is no column reversing, so we have the result. If
$t=0$, then the number of elements of $E_2$ is $\sharp
E_1-\sharp F$, where $F=E_1\cap\{\bullet\mbox{-nw-diagrams of } \mathcal S^{ab}_I
\}.$ By proposition \ref{condabelienD}, we have $F=Q \cup \widetilde
P$, so by the consideration
above, we have $\sharp F=2\sharp Q$. It follows that $\sharp
E_1+\sharp E_2=2\sharp P$.

For the last case $t=2$, the shape of $T_{D_l}^I$ is
$T_{l-r,l-r}(l_1,\dots, l_{s-1})$. If $S\in P$, then $S$ is 
included in the left hand half of
$T_{D_l}^{I}$, so by proposition \ref{compteab}, we have 
$$\sharp P=2^{l-r-1}+\sum_{j=1}^{s-1} 
\left(\begin{array}{c} l-r-1 \\l_j-1\end{array}\right).
$$ 
We have $E_2=\emptyset$, and $Q_i$ is defined for $i=l-r,\dots,
2(l-r)-1$. Set $P_i=\{S\in P;
\tau_1^S=2(l-r)-1-i\}$. We then have : 
$$
\displaystyle
Q=\bigcup_{i=l-r}^{2(l-r)-1} Q_i \ \mbox{ and }\ 
P=\bigcup_{i=l-r}^{2(l-r)-1} P_i.
$$
As above, for $i=l-r,\dots,2(l-r)-1$, we have an obvious bijection between $P_i$ and
$Q_i$ given by the adding or removing of boxes $(1,2(l-r)-i),
\dots,(1,i)$. Therefore $\sharp P=\sharp Q$ and the result follows.
\end{proof}

\begin{remarque}
All the results above depend on the numbering of simple roots.
\end{remarque}

%----------------------------------------------------------------------------------------------------------------------
\section{Remarks}

\subsection{Exceptional types}
In the exceptional types $E$, $F$ and $G$, the number of ad-nilpotent and
abelian ideals has been determined by using GAP 4. 

The following tables give the cardinality of $\calf_I$ and $\mathcal{A}b_{I}$ for the types $F_4$ and $G_2$. The subset $I$ of $\Pi$ is described by the symbol $\bullet$ in the Dynkin diagram without arrow.

%% Table F4

$$
\begin{array}{||c|c|c||c|c|c||}
\hline
I & \sharp\calf_{I} & \sharp \mathcal{A}b_{I}  & 
I &  \sharp\calf_{I} &  \sharp\mathcal{A}b_{I} \\
\hline
\begin{smallmatrix}& & 
\circ &\circ & \circ &\circ & \\
\end{smallmatrix}
&105&16& 
\begin{smallmatrix}& & \\
\bullet &\circ &\circ &\circ &\\
\end{smallmatrix}
&24&6\\
\hline
\begin{smallmatrix}& & 
\circ &\bullet & \circ &\circ &\\
\end{smallmatrix}
&35&12& 
\begin{smallmatrix}& & \\
\circ &\circ & \bullet &\circ & \\
\end{smallmatrix}
&32&10 \\
\hline
\begin{smallmatrix}& & 
\circ &\circ &\circ &\bullet &\\
\end{smallmatrix}
&49&9&
\begin{smallmatrix}& & \\
\bullet &\bullet & \circ &\circ & \\
\end{smallmatrix}
&10&5 \\
\hline
\begin{smallmatrix}& & 
\bullet &\circ &\bullet &\circ &\\
\end{smallmatrix}
&8&4&
\begin{smallmatrix}& & \\
\bullet &\circ & \circ &\bullet &\\
\end{smallmatrix}
&12&4 \\
\hline
\begin{smallmatrix}& & 
\circ &\bullet &\bullet &\circ &\\
\end{smallmatrix}
&14&7 &
\begin{smallmatrix}& & \\
\circ &\bullet & \circ &\bullet & \\
\end{smallmatrix}
&14&6 \\ 
\hline 
\begin{smallmatrix}& & 
\circ &\circ &\bullet &\bullet &\\
\end{smallmatrix}
&10&4&
\begin{smallmatrix}& & \\
\bullet &\bullet & \bullet &\circ & \\
\end{smallmatrix}
&4&3 \\
\hline
\begin{smallmatrix}& & 
\bullet &\bullet &\circ &\bullet &\\
\end{smallmatrix}
&5&3 &
\begin{smallmatrix}& & \\
\bullet &\circ & \bullet &\bullet &\\
\end{smallmatrix}
&3&2 \\
\hline
\begin{smallmatrix}& & 
\circ &\bullet &\bullet &\bullet &\\
\end{smallmatrix}
&3&2&
\begin{smallmatrix}& & \\
\bullet &\bullet & \bullet &\bullet & \\
\end{smallmatrix}
&1&1 \\
\hline
\end{array}
$$
$$
$$
where we use the following orientation for the Dynkin diagram of $F_4$ :
$$
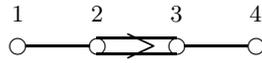
\begin{Dynkin}
\Dbloc{\Dcirc\Deast\Dtext{t}{1}}
\Dbloc{\Dcirc\Dwest\Ddoubleeast\Dtext{t}{2}}
\Drightarrow
\Dbloc{\Dcirc\Ddoublewest\Deast\Dtext{t}{3}}
\Dbloc{\Dcirc\Dwest\Dtext{t}{4}}
\end{Dynkin}
$$

%% Table G2

$$
\begin{array}{||c|c|c||c|c|c||}
\hline
I & \sharp\calf_{I} & \sharp \mathcal{A}b_{I}  & 
I &  \sharp\calf_{I} &  \sharp\mathcal{A}b_{I} \\
\hline
\begin{smallmatrix}& & 
\circ &\circ & \\
\end{smallmatrix}
&8&4& 
\begin{smallmatrix}& & \\
\bullet &\circ &\\
\end{smallmatrix}
&3&2\\
\hline
\begin{smallmatrix}& & 
\circ &\bullet & \\
\end{smallmatrix}
&4&3& 
\begin{smallmatrix}& & \\
\bullet &\bullet &\\
\end{smallmatrix}
&1&1\\
\hline
\end{array}
$$
where we use the following orientation for the Dynkin diagram of $G_2$ :

$$
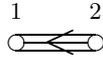
\begin{Dynkin}
\Dbloc{\Dcirc\Deast\Ddoubleeast\Dtext{t}{1}}
\Dleftarrow
\Dbloc{\Dcirc\Dwest\Ddoublewest\Dtext{t}{2}}
\end{Dynkin}
$$

\subsection{Relation with antichains}
For $\Phi\in \calf_{\emptyset}$, let 
$$
\Phi_{min}=\{\beta\in\Phi;
\beta-\alpha\not\in \Phi, \mbox{ for all }\alpha\in \Delta^+\}
$$
be the set of minimal roots of $\Phi$, also called an antichain of 
$(\Delta^+,\leqslant)$, see \cite{P}. It is clear that each antichain
corresponds to an element of $\calf_{\emptyset}$ and vice versa.

By a similar proof as in \cite{P}, we obtain the following proposition: 

\begin{prop}
Let $I \subset \Pi$ be of cardinality $r$ and $\Phi\in \calf_I$, then we have $\sharp \Phi_{min} \leqslant l-r$. 
\end{prop}

\begin{proof}
Let $I \subset \Pi$ be of cardinality $r$ and $\Phi\in \calf_I$. Set $\Gamma=\Phi_{min}\cup I=\{\gamma_1,\dots,\gamma_t\}.$ Let $\gamma_i,\gamma_j\in\Gamma$, then $\gamma_i-\gamma_j\not\in\Delta$ by the definition of $\Phi_{min}$ and the fact that $\Phi\in\calf_I$. Thus the angle between any pair of distinct elements of $\Gamma$ is non acute and since all the $\gamma_i$'s lie in an open half space of $V$, they are linearly independent. Consequently, we have $\sharp \Gamma\leqslant r$, and hence $\sharp \Phi_{min}\leqslant l-r$.
\end{proof}

\begin{remarques}
{\rm{(i)}} Recall from \cite{CP}, that an antichain $\Gamma\subset\Delta^+$ is of cardinality $l$ if and only if $\Gamma=\Pi$. This result has no equivalence in the general parabolic case. For example, in $B_2$, the set $\Phi=\{\alpha_1+2\alpha_2\}$ is an ad-nilpotent ideal of $\plie_{\alpha_1}$ such that $\Phi_{min}=\Phi$ and $\sharp\Phi_{min}=1$. 

{\rm{(ii)}} Let $\glie$ be of type $A_l$. Let $I \subset \Pi$ be of cardinality $r$ and $\Phi\in \calf_I$, then it is possible to show that $\sharp \Phi_{min}=l-r$ if and only if $\Phi_{min}=\Pi\setminus I$.

\end{remarques}

\end{document}